\RequirePackage{fix-cm}
\documentclass[smallextended]{svjour3}       
\smartqed  

\usepackage{graphicx}
\usepackage{amsmath,amsthm,amsfonts,amscd,amssymb,bm,mathrsfs}
\usepackage{enumerate}
\usepackage{graphicx,xcolor,color}
\usepackage{epstopdf}
\usepackage{subfigure}
\usepackage{multirow}
\usepackage{cases}
\DeclareMathOperator*{\argmin}{arg\,min}

\newcommand{\B}{{\mathbb{B}}}
\newcommand{\M}{{\mathbb{M}}}
\newcommand{\K}{{\mathbb{K}}}
\newcommand{\C}{{\mathbb{C}}}
\newcommand{\R}{{\mathbb{R}}}
\newcommand{\Q}{{\mathbb{Q}}}
\newcommand{\Pro}{{\mathbb{P}}}
\newcommand{\E}{{\mathbb{E}}}
\newcommand{\Su}{{\mathbb{S}}}
\newcommand{\Aa}{{\mathbb{A}}}
\newcommand{\Bb}{{\mathbb{B}}}

\newcommand{\I}{{\mathbb{I}}}
\newcommand{\N}{{\mathcal{N}}}
\newcommand{\loss}{{\mathcal{L}}}
\newcommand{\obj}{\Psi}

\newcommand{\normm}[1]{{\left\vert\kern-0.25ex\left\vert\kern-0.25ex\left\vert #1
\right\vert\kern-0.25ex\right\vert\kern-0.25ex\right\vert}}
\newcommand{\inm}[2]{\langle\langle #1,#2 \rangle\rangle}

\newtheorem{Lemma}{Lemma}
\newtheorem{Definition}{Definition}
\newtheorem{Proposition}{Proposition}
\newtheorem{Remark}{Remark}
\newtheorem{Theorem}{Theorem}
\newtheorem{Corollary}{Corollary}

\begin{document}

\title{Low-rank matrix estimation via nonconvex optimization methods in multi-response errors-in-variables regression
}

\titlerunning{Low-rank matrix estimation via nonconvex optimization methods}        

\author{Xin Li        \and
        Dongya Wu 
}

\authorrunning{Xin Li and Dongya Wu} 

\institute{X. Li \at
             School of Mathematics, Northwest University, Xi’an, 710069, P. R. China \\
              \email{lixin@nwu.edu.cn}
          \and
          D. Wu \at
              School of Information Science and Technology, Northwest University, Xi’an, 710069, P. R. China \\
              \email{wudongya@nwu.edu.cn}
}

\date{Received: date / Accepted: date}

\maketitle

\begin{abstract}
Noisy and missing data cannot be avoided in real application, such as bioinformatics, economics and remote sensing. Existing methods mainly focus on  linear errors-in-variables regression, while relatively little attention is paid for the case of multivariate responses, and how to achieve consistent estimation under corrupted covariates is still an open question. In this paper, a nonconvex error-corrected estimator is proposed for the matrix estimation problem in the multi-response errors-in-variables regression model. Statistical and computational properties for global solutions of the estimator are analysed. In the statistical aspect, the nonasymptotic recovery bound for all global solutions of the nonconvex estimator is established. In the computational aspect, the proximal gradient method is applied to solve the nonconvex optimization problem and proved to linearly converge to a near-global solution. Sufficient conditions are verified in order to obtain probabilistic consequences for specific types of measurement errors by virtue of random matrix analysis. Finally, simulation results demonstrate the theoretical properties and show nice consistency under high-dimensional scaling.
\end{abstract}

\section{Introduction}\label{sec-intro}

Massive data sets have posed a variety of challenges to the field of statistics and machine learning in recent decades. In view of these challenges, researchers have developed different classes of statistical models and numerical algorithms to deal with the complexities of modern data; see the books \cite{buhlmann2011statistics,wainwright2019high} for an overall review.

In standard formulations of statistical inference problems, it is assumed that the collected data are clean enough. However, this hypothesis is neither realistic nor reasonable, since in many real-world problems, due to instrumental constraints or a lack of observation, the collected data may always be corrupted and tend to be noisy or partially missing. Hence, measurement errors cannot be avoided in general. There have been a variety of researches focusing on models with corrupted data for regression problems in low dimensions; see, e.g., \cite{bickel1987efficient,carroll2006measurement} and references therein. As for the high-dimensional scenario, it has been pointed in \cite{sorensen2015measurement} that misleading inference results can still only be obtained if the method for clean data is applied naively to the corrupted data. Thus it is necessary to take measurement errors into consideration and develop new error-corrected methods for high-dimensional models.

Recently, some regularization methods have been proposed to deal with high-dimensional errors-in-variables regression. For example, Loh and Wainwright proposed a nonconvex Lasso-type estimator via substituting the unobserved matrices involved in the ordinary least squares loss function with unbiased surrogates, and established statistical errors for global solutions \cite{loh2012high,loh2015regularized}. 
In order to avoid the nonconvexity, Datta and Zou developed the convex conditional Lasso (CoCoLasso) by defining the nearest positive semi-definite matrix \cite{datta2017cocolasso}. CoCoLasso thus enjoys the benefits of convex optimization and is shown to possess nice estimation accuracy in linear regression. Rosenbaum and Tsybakov proposed a modified form of the Dantzig selector \cite{candes2007the}, called matrix uncertainty selector (MUS) \cite{rosenbaum2010sparse}. 
Further development of MUS included modification to achieve statistical consistency, and generalization to deal with the cases of unbounded and dependent measurement errors as well as generalized linear models \cite{belloni2016an,belloni2017linear,rosenbaum2013improved,sorensen2018covariate}. Li et al. investigated a general nonconvex regularized $M$-estimator, which can be applied to errors-in-variables sparse linear regression, and analysed the statistical and computational properties \cite{li2020sparse}. Li et al. studied statistical inference problems and proposed a corrected decorrelated score test and a score type estimator \cite{li2021inference}. 

However, until now, relatively little attention is paid for the case of multivariate responses. Though a simple and natural idea is to vectorize both the response matrix and the coefficient matrix so that methods for the univariate response case can be directly applied, it may ignore the particular low-dimensional structures of the coefficient matrix such as low-rank and row/column sparsity as well as the multivariate nature of the responses \cite{izenman2008modern}. Moreover, the multi-response linear regression model has a substantial wider application than that of the univariate case in modern large-scale association analysis, such as genome-wide association studies (GWASs) \cite{li2019reliable} and social network analyses \cite{mcgillivray2020estimating}. 

In this work, we shall deal with multi-response errors-in-variables regression under high-dimensional scaling. More precisely, we consider the multi-response regression model $Y=X\Theta^*+\epsilon$, where $\Theta^*\in \R^{d_1\times d_2}$ is the unknown underlying matrix parameter. The covariate matrix $X\in \R^{n\times d_1}$ does not need to be fully observed, instead, two types of errors will be discussed, including the additive noise and missing data cases. A nonconvex error-corrected estimator is constructed with its statistical and computational properties established.

As has been pointed out in \cite{candes2010thepower} that low-rank matrix regression is intrinsically different from sparse linear regression in that regularity conditions to ensure recovery are not satisfied in this case. Therefore, conditions used in linear regression cannot be used to analyse errors-in-variables matrix regression. More general conditions called restricted strong convexity/smoothness (RSC/RSM) have been introduced in \cite{negahban2011estimation} to ensure statistical and computational properties of low-rank matrix regression without measurement error.

To date, however, an open question is whether or not an appropriate form of the RSC/RSM conditions holds for errors-in-variables matrix regression, and this leads to the first challenge of this article. The second challenge is that to what extent do the plug-in surrogates approximate the unobserved matrices. If the approximation degree is high, then together with suitable regularity conditions, it would be possible to derive nonasymptotic recovery bounds and to develop efficient algorithms for errors-in-variables matrix regression models.

The major contributions of this article are threefold. First, it is proved that with overwhelming probability, a specialized form of the RSC/RSM conditions holds for errors-in-variables matrix regression, and the surrogates well approximate the unobserved matrices by virtue of a suitable deviation control condition. Results and proofs corresponding to specific error models make use of nontrivial matrix analysis on concentration inequalities and random matrix theory (cf. Propositions \ref{prop-add-rs}--\ref{prop-mis-devia}), which cannot be easily generalized from linear regression. Second, recovery bounds are provided for global solutions of the proposed nonconvex error-corrected estimator to imply statistical consistency (cf. Theorem \ref{thm-stat}). Last but not least, the proximal gradient method \cite{nesterov2007gradient} is applied to solve the optimization problem and is proved to converge linearly to a local solution which lies in a small neighborhood of all global optima. This local solution essentially performs as well as any global solution in the sense of statistical precision and thus serves as a nice approximation to any global solution (cf. Theorem \ref{thm-algo}).

The remainder of this article is organized as follows. In Sect. \ref{sec-prob}, we provide background on the high-dimensional multi-response errors-in-variables regression model and then construct an error-corrected estimator. Some regularity conditions are imposed to facilitate the analysis. In Sect. \ref{sec-main}, we establish our main results on statistical recovery bounds and computational convergence rates. In Sect. \ref{sec-conse}, probabilistic consequences for specific errors-in-variables models are obtained by verifying the regularity conditions. In Sect. \ref{sec-simul}, several numerical experiments are performed to demonstrate theoretical properties. Conclusions and future work are discussed in Sect. \ref{sec-concl}. Technical proofs are deferred to Appendix. 

While we are working on this paper, we became aware of
an independent related work \cite{wu2020scalable}. Our contributions are substantially different from theirs, in that our theory consists of statistical recovery   bounds and algorithmic convergence rates, whereas they only provide recovery bounds. In addition, our theory is applicable to both the exact low-rank case and the near low-rank case, while they only consider the exact low-rank case. Finally, we provide rigorous proofs showing that the regularity conditions are satisfied with high probability for different types of errors, whereas the applicability of their technical conditions was not established.

We end this section by introducing useful notations. For $d\geq 1$, let $\I_d$ stand for the $d\times d$ identity matrix. For a matrix $X\in \R^{n\times d}$, let $X_{ij}\ (i=1,\dots,n,j=1,2,\cdots,d)$ denote its $ij$-th entry, $X_{i\cdot}\ (i=1,\dots,n)$ denote its $i$-th row, and $X_{\cdot j}\ (j=1,2,\cdots,d)$ denote its $j$-th column. When $X$ is a square matrix, i.e., $n=d$, we use diag$(X)$ stand for the diagonal matrix with its diagonal elements equal to $X_{11},X_{22},\cdots,X_{dd}$. We write $\lambda_{\text{min}}(X)$ and $\lambda_{\text{max}}(X)$ to denote the minimal and maximum eigenvalues of a matrix $X$, respectively. For a matrix $\Theta\in \R^{d_1\times d_2}$, define $d=\min\{d_1,d_2\}$, and denote its singular values in decreasing order by $\sigma_1(\Theta)\geq \sigma_2(\Theta)\geq \cdots \sigma_d(\Theta)\geq 0$. We use $\normm{\cdot}$ to denote different types of matrix norms based on singular values, including the nuclear norm $\normm{\Theta}_*=\sum_{j=1}^{d}\sigma_j(\Theta)$, the spectral or operator norm $\normm{\Theta}_{\text{op}}=\sigma_1(\Theta)$, and the Frobenius norm $\normm{\Theta}_F=\sqrt{\text{trace}(\Theta^\top\Theta)}=\sqrt{\sum_{j=1}^{d}\sigma_j^2(\Theta)}$. All vectors are column vectors following classical mathematical convention. For a pair of matrices $\Theta$ and $\Gamma$ with equal dimensions, we let $\inm{\Theta}{\Gamma}=\text{trace}(\Theta^\top \Gamma)$ denote the trace inner product on matrix space. For a function $f:\R^d\to \R$, $\nabla f$ is used to denote a gradient or subgradient depending on whether $f$ is differentiable or nondifferentiable but convex, respectively.

\section{Problem setup}\label{sec-prob}

In this section, we begin with background on the high-dimensional multi-response errors-in-variables regression, and then a precise description on the proposed nonconvex error-corrected estimation method. Finally, we introduce some regularity conditions that play key roles in the following analysis.

\subsection{Model setting}

Consider the high-dimensional multi-response regression model which links the response vector $Y_{i\cdot}\in \R^{d_2}$ to a covariate vector $X_{i\cdot}\in \R^{d_1}$
\begin{equation}\label{eq-model}
Y_{i\cdot}={\Theta^*}^\top X_{i\cdot}+\epsilon_{i\cdot}\qquad \text{for}\ i=1,2,\cdots,n,
\end{equation}
where $\Theta^*\in \R^{d_1\times d_2}$ is the unknown parameter matrix, and $\epsilon_{i\cdot}\in \R^{d_2}$ is the observation noise independent of $X_{j\cdot}$ ($\forall\ i,j$).
Model \eqref{eq-model} can be written in a more compact form using matrix notation. Particularly, define the multi-response matrix $Y=(Y_{1\cdot},Y_{2\cdot},\cdots,Y_{n\cdot})^\top\in \R^{n\times d_2}$ with similar definitions for the covariate matrix $X\in \R^{n\times d_1}$ and the noise matrix $\epsilon\in \R^{n\times d_2}$ in terms of $\{X_{i\cdot}\}_{i=1}^n$ and $\{\epsilon_{i\cdot}\}_{i=1}^n$, respectively. Then model \eqref{eq-model} is re-written as
\begin{equation}\label{eq-model-matrix}
Y=X\Theta^*+\epsilon.
\end{equation}

This work is applicable to the high-dimensional case where the number of covariates or responses (i.e., $d_1$ or $d_2$) may be possibly more than the sample size $n$. It is already known that consistent estimation cannot be obtained under this high-dimensional regime unless the model is imposed with additional structures, such as low-rank in matrix estimation problems.
In the following, unless otherwise specified, $\Theta^*\in \R^{d_1\times d_2}$ is assumed to be of either exact low-rank, i.e., it has rank far less than $\min\{d_1,d_2\}$, or near low-rank, meaning that it can be approximated by a matrix of low-rank perfectly. One popular way to measure the degree of low-rank is to use the matrix $\ell_q$-ball (accurately speaking, when $q\in[0,1)$, these sets are not real ``ball''s due to the nonconvexity), which is defined as, for $q\in [0,1]$, and a radius $R_q>0$,
\begin{equation}\label{eq-lqball}
\B_q(R_q):=\{\Theta\in \R^{d_1\times d_2}\big|\sum_{i=1}^{d}|\sigma_i(\Theta)|^q\leq R_q\},
\end{equation}
where $d=\min\{d_1,d_2\}$. Note that the matrix $\ell_0$-ball corresponds to the case of exact low-rank, meaning that the rank of a matrix is at most $R_0$; while the matrix $\ell_q$-ball for $q\in (0,1]$ corresponds to the case of near low-rank, which enforces a certain decay rate on the ordered singular values of the matrix $\Theta\in \B_q(R_q)$. In this paper, we fix $q\in [0,1]$, and assume that the true parameter $\Theta^*\in \B_q(R_q)$ unless otherwise specified.

In standard formulations, the true covariates $X_{i\cdot}$'s are assumed to be fully-observed. However, this assumption is not realistic for many real-world applications, in which the covariates may be observed with measurement errors. We use $Z_{i\cdot}$'s to denote noisy versions of the corresponding $X_{i\cdot}$'s. It is usually assumed that $Z_{i\cdot}$ is linked to $X_{i\cdot}$ via some conditional distribution as follows
\begin{equation}\label{eq-link}
Z_{i\cdot}\sim \Q(\cdot|X_{i\cdot})\quad \text{for}\ i=1,2,\cdots,n.
\end{equation}
As has been discussed in \cite{loh2012high,loh2015regularized}, there are mainly two types of errors:
\begin{enumerate}[(a)]
\item Additive noise: For each $i=1,2,\cdots,n$, one observes $Z_{i\cdot}=X_{i\cdot}+W_{i\cdot}$, where $W_{i\cdot}\in \R^{d_1}$ is a random vector independent of $X_{i\cdot}$ with mean 0 and known covariance matrix $\Sigma_w$.
\item Missing data: For each $i=1,2,\cdots,n$, one observes a random vector $Z_{i\cdot}\in \R^{d_1}$, such that for each $j=1,2\cdots,d_1$, we independently observe $Z_{ij}=X_{ij}$ with probability $1-\rho$, and $Z_{ij}=*$ with probability $\rho$, where $\rho\in [0,1)$.
\end{enumerate}
Throughout this paper, unless otherwise specified, we assume a Gaussian model on the covariate and error matrices. Specifically, the matrices $X$, $W$ and $\epsilon$ are assumed to be random matrices with independent and identically distributed (i.i.d.) rows sampled from Gaussian distributions $\N(0,\Sigma_x)$, $\N(0,\sigma_w^2\mathbb{I}_{d_1})$ and $\N(0,\sigma_\epsilon^2\mathbb{I}_{d_2})$, respectively. Then it follows that $\Sigma_w=\sigma_w^2\mathbb{I}_{d_1}$. When the noise covariance $\Sigma_w$ is unknown, one needs to estimate it by additional experiments; see, e.g., \cite{carroll2006measurement}. For example, a simple method to estimate $\Sigma_w$ is by independent observations of the noise blankly. Specifically, suppose a matrix $W_0\in \R^{n\times d}$ is independently observed with $n$ i.i.d. vectors of noise. Then the population covariance $\frac1n W_0^\top W_0$ is used as the estimate of $\Sigma_w$. Some other variant of this method in are also discussed in \cite{carroll2006measurement}. Define
$Z=(Z_{1\cdot},Z_{2\cdot},\cdots,Z_{n\cdot})^\top$. Thus $Z$ is the observed covariate matrix, which involves certain types of measurement errors, and will be specified according to the context.

\subsection{Error-corrected M-estimators}\label{sec-esti}

When the covariate matrix $X$ is correctly obtained, previous literatures have proposed various methods for rank-constrained problems; see, e.g., \cite{negahban2011estimation,recht2010guaranteed,zhou2014regularized} and references therein. Most of the estimators are formulated as solutions to certain semidefinite programs (SDPs) based on the nuclear norm regularization. Recall that for a matrix $\Theta\in \R^{d_1\times d_2}$, the nuclear or trace norm is defined by $\normm{\Theta}_*=\sum_{j=1}^{d}\sigma_j(\Theta)\ (d=\min\{d_1,d_2\})$, corresponding to the sum of its singular values. For instance, given model \eqref{eq-model-matrix}, define $N=nd_2$, a commonly-used estimator is based on solving the following SDP:
\begin{equation}\label{eq-wain}
\hat{\Theta}\in \argmin_{\Theta\in \Su}\left\{\frac{1}{2N}\normm{Y-X\Theta}_F^2+\lambda_N\|\Theta\|_*\right\},
\end{equation}
where $\Su$ is some subset of $\R^{d_1\times d_2}$, and $\lambda_N>0$ is a regularization parameter. The nuclear norm of a matrix offers a natural convex relaxation of the rank constraint, analogous to the $\ell_1$ norm as a convex relaxation of the cardinality for a vector. The statistical and computational properties of the nuclear norm regularization methods have been studied deeply and applied widely in various fields, such as matrix completion, matrix decomposition and so on \cite{agarwal2012noisy,negahban2012restricted}. Other regularizers including the elastic net, SCAD and MCP, which were proposed to deal with sparse linear regression, have also been adapted to low-rank matrix estimation problems in terms of matrix singular value vectors; see \cite{zhou2014regularized} for a detailed investigation.

Note that \eqref{eq-wain} can be re-written as follows
\begin{equation}\label{eq-rewain}
\hat{\Theta}\in \argmin_{\Theta\in \Su}\left\{\frac{1}{2N}\inm{X^\top X\Theta}{\Theta}-\frac{1}{N}\inm{X^\top Y}{\Theta}+\lambda_N\normm{\Theta}_*\right\}.
\end{equation}
Recall the relation $N=nd_2$, the SDP optimization problem \eqref{eq-rewain} is then transformed to
\begin{equation}\label{eq-rewain1}
\hat{\Theta}\in \argmin_{\Theta\in \Su}\frac{1}{d_2}\left\{\frac{1}{2n}\inm{X^\top X\Theta}{\Theta}-\frac{1}{n}\inm{X^\top Y}{\Theta}+\lambda_n\normm{\Theta}_*\right\},
\end{equation}
where we have defined $\lambda_n=d_2\lambda_N$.
In the errors-in-variables cases, the quantities $\frac{X^\top X}{n}$
and $\frac{X^\top Y}{n}$ in \eqref{eq-rewain1} are both unknown, meaning that this estimator does not work any more. However, this transformation still provides some useful intuition for the construction of estimation methods by virtue of the plug-in principle proposed in \cite{loh2012high}. Specifically, given a set of samples, one way is to find suitable estimates of the quantities $\frac{X^\top X}{n}$ and $\frac{X^\top Y}{n}$ that are adapted to the cases of additive noise and/or missing data. 

Let $(\hat{\Gamma},\hat{\Upsilon})$ denote estimates of  $(\frac{X^\top X}{n},\frac{X^\top Y}{n})$. Inspired by \eqref{eq-rewain1} and ignoring the constant, we propose the following estimator to solve the low-rank estimation problem in errors-in-variables regression:
\begin{equation}\label{eq-mine}
\hat{\Theta}\in \argmin_{\Theta\in \Su}\left\{\frac{1}{2}\inm{\hat{\Gamma}\Theta}{\Theta}-\inm{\hat{\Upsilon}}{\Theta}+\lambda\normm{\Theta}_*\right\},
\end{equation}
where $\lambda>0$ is the regularization parameter. The feasible region is specialized as $\Su=\{\Theta\in \R^{d_1\times d_2}\big|\normm{\Theta}_*\leq \omega\}$, and the parameter $\omega>0$ must be chosen carefully to guarantee $\Theta^*\in \Su$. We include this side constraint is because of the nonconvexity of the estimator, which will be explained in detail in the last paragraph of this section. Then any matrix $\Theta\in \Su$ will also satisfy $\normm{\Theta}_*\leq \omega$, and since the optimization object function is continuous, it is guaranteed by Weierstrass extreme value theorem that a global solution $\hat{\Theta}$ always exists. Though when the covariates are clean, either the regularization term or the side constraint is enough to impose low-rank constraints, in errors-in-variables regression models, these two terms play completely different roles in \eqref{eq-mine}. Specifically, the regularization term is to impose low-rank constraints, while the side constraint is used to ensure the existence of the global solution.

Note that the estimator \eqref{eq-mine} is just a general expression, and concrete expressions for $(\hat{\Gamma},\hat{\Upsilon})$ still need to be determined. For the specific additive noise and missing data cases, as discussed in \cite{loh2012high}, an unbiased choice of the pair $(\hat{\Gamma},\hat{\Upsilon})$ is given respectively by
\begin{align}
\hat{\Gamma}_{\text{add}}&:=\frac{Z^\top Z}{n}-\Sigma_w \quad \mbox{and}\quad \hat{\Upsilon}_{\text{add}}:=\frac{Z^\top Y}{n}, \label{sur-add}\\
\hat{\Gamma}_{\text{mis}}&:=\frac{\tilde{Z}^\top \tilde{Z}}{n}-\rho\cdot\text{diag}\left(\frac{\tilde{Z}^\top \tilde{Z}}{n}\right) \quad \mbox{and}\quad \hat{\Upsilon}_{\text{mis}}:=\frac{\tilde{Z}^\top Y}{n}\quad \left(\tilde{Z}=\frac{Z}{1-\rho}\right). \label{sur-mis}
\end{align}

Under the high-dimensional regime $(n\ll d_1)$, the surrogate matrices $\hat{\Gamma}_{\text{add}}$ and $\hat{\Gamma}_{\text{mis}}$ in \eqref{sur-add} and \eqref{sur-mis} are always negative definite; indeed, both the matrices $Z^{\top}Z$ and $\tilde{Z}^{\top} \tilde{Z}$ are with rank at most $n$, and then the positive definite matrices $\Sigma_{w}$ and $\rho\cdot \text{diag}\left(\frac{\tilde{Z}^\top \tilde{Z}}{n}\right)$ are subtracted to arrive at the estimates $\hat{\Gamma}_{\text{add}}$ and $\hat{\Gamma}_{\text{mis}}$, respectively. Therefore, the above estimator \eqref{eq-mine} involves solving a nonconvex optimization problem. Due to the nonconvexity, it is generally impossible to obtain a global solution through a polynomial-time algorithm. Nevertheless, this issue is not significant in our setting, and we shall establish that a simple proximal gradient algorithm converges linearly to a matrix lying extremely close to any global optimum of the problem \eqref{eq-mine} with high probability.

\subsection{Regularity conditions}

Now we impose some regularity conditions on the surrogate matrices $\hat{\Gamma}$ and $\hat{\Upsilon}$, which will be beneficial to the statistical and computational analysis for the nonconvex estimator \eqref{eq-mine}.

Regularity conditions called RSC/RSM have been adopted to analyse theoretical properties of low-rank matrix regression, and is applicable when the loss function is nonquadratic or nonconvex; see, e.g., \cite{li2020sparse,loh2015regularized,negahban2011estimation}. For linear/matrix regression without measurement errors, it has been shown that the RSC/RSM conditions are satisfied by various types of random matrices with high probability \cite{agarwal2012fast,negahban2011estimation}. 


However, it is still unknown whether or not an appropriate form of RSC/RSM holds for errors-in-variables matrix regression. In this paper, we provide a positive answer for this question, and propose the following general RSC/RSM conditions. Verification for specific types of measurement error serves as the first challenge of this article, which involves probabilistic arguments under high-dimensional scaling and will be established in Section \ref{sec-conse}.

\begin{Definition}[Restricted strong convexity]
The matrix $\hat{\Gamma}$ is said to satisfy a restricted strong convexity with curvature $\alpha_1>0$ and tolerance $\tau(n,d_1,d_2)>0$ if
\begin{equation}\label{eq-rsc}
\inm{\hat{\Gamma}\Delta}{\Delta}\geq \alpha_1\normm{\Delta}_F^2-\tau(n,d_1,d_2)\normm{\Delta}_*^2\quad \emph{for all}\ \Delta\in R^{d_1\times d_2}.
\end{equation}
\end{Definition}

\begin{Definition}[Restricted strong smoothness]
The matrix $\hat{\Gamma}$ is said to satisfy a restricted strong smoothness with smoothness $\alpha_2>0$ and tolerance $\tau(n,d_1,d_2)>0$ if
\begin{equation}\label{eq-rsm}
\inm{\hat{\Gamma}\Delta}{\Delta}\leq \alpha_2\normm{\Delta}_F^2+\tau(n,d_1,d_2)\normm{\Delta}_*^2\quad \emph{for all}\ \Delta\in R^{d_1\times d_2}.
\end{equation}
\end{Definition}

In the following, we shall adopt the shorthand $\tau$ for $\tau(n,d_1,d_2)$ to ease the notation. 


Recall that $(\hat{\Gamma},\hat{\Upsilon})$ are estimates for the unknown quantities $(\frac{X^\top X}{n},\frac{X^\top Y}{n})$. Then it is natural to ask that to what extent the estimates approximate the true unknown quantities. \cite{loh2012high} has shown that in linear errors-in-variables regression, the approximation is nice under high-dimensional scaling. However, this question has not been settled down for multi-response errors-in-variables regression. In the following, we show that a suitable deviation bound is enough to control the approximation degree.
Specifically, assume that there is some function $\phi(\Q,\sigma_\epsilon)$, hinging on the two sources of errors in our setting: the conditional distribution $\Q$ (cf. \eqref{eq-link}) which links the true covariates $X_{i\cdot}$
to the corrupted versions $Z_{i\cdot}$ and the standard deviation $\sigma_\epsilon$ of the observation noise $\epsilon_{i\cdot}$. With this notation, the following deviation condition is required:
\begin{equation}\label{eq-devia}
\normm{\hat{\Upsilon}-\hat{\Gamma}\Theta^*}_{\text{op}}\leq \phi(\Q,\sigma_\epsilon)\sqrt{\frac{\max (d_1,d_2)}{n}}.
\end{equation}
Note that when there exists no measurement error, and thus $\hat{\Gamma}=\frac{X^\top X}{n}$, $\hat{\Upsilon}=\frac{X^\top Y}{n}$ hold trivially, the left hand of \eqref{eq-devia} turns to $\frac{1}{n}\normm{X^\top\epsilon}_{\text{op}}$—a quantity which has been used to determine the regularization parameter in matrix estimation \cite{negahban2011estimation}. We shall see in the following that \eqref{eq-devia} also plays a part in deciding the regularization parameter in the errors-in-variables model. Though the idea of derivation condition \eqref{eq-devia} is similar to that of \cite{loh2012high}, the proof is intrinsically different from that of linear regression due to matrix structure of the underlying parameter, and will be established for various forms of errors in Section \ref{sec-conse}.


\section{Main results}\label{sec-main}

In this section, we establish our main results including statistical guarantee on the recovery bound for the nonconvex estimator \eqref{eq-mine} and computational guarantee on convergence rates for the proximal gradient algorithm. These results are deterministic in nature. Consequences for the concrete additive noise and missing data cases will be discussed in the next section. 


Before we proceed, some additional notations are required to facilitate the analysis of exact/near low-rank matrices. First let $\obj(\Theta)=\frac{1}{2}\inm{\hat{\Gamma}\Theta}{\Theta}-\inm{\hat{\Upsilon}}{\Theta}+\lambda\normm{\Theta}_*$ denote the objective function to be minimized, and $\loss(\Theta)=\frac{1}{2}\inm{\hat{\Gamma}\Theta}{\Theta}-\inm{\hat{\Upsilon}}{\Theta}$ denote the loss function. Then it follows that $\obj(\Theta)=\loss(\Theta)+\lambda\normm{\Theta}_*$.

Note that the parameter matrix $\Theta^*$ has a singular value
decomposition of the form $\Theta^*=UDV^\top$, where $U\in \R^{d_1\times d}$ and $V\in \R^{d_2\times d}$
are orthonormal matrices with $d=\min\{d_1,d_2\}$ and without loss of generality, assume that $D$ is diagonal with singular values in nonincreasing order, i.e., $\sigma_1(\Theta^*)\geq \sigma_2(\Theta^*)\geq \cdots \sigma_d(\Theta^*)\geq 0$. For each integer $r\in \{1, 2,\cdots,d\}$, let $U^r\in \R^{d_1\times r}$ and $V^r\in \R^{d_2\times r}$ be the sub-matrices consisting of singular vectors corresponding to the first $r$ largest
singular values of $\Theta^*$. Then two subspaces of $\R^{d_1\times d_2}$ associated with $\Theta^*$ are defined as follows:
\begin{subequations}\label{eq-sub}
\begin{align}\
\Aa(U^r,V^r)&:=\{\Delta\in \R^{d_1\times d_2}\big|\text{row}(\Delta)\subseteq \text{col}(V^r), \text{col}(\Delta)\subseteq \text{col}(U^r)\}\quad \text{and}\label{eq-sub1}\\
\Bb(U^r,V^r)&:=\{\Delta\in \R^{d_1\times d_2}\big|\text{row}(\Delta)\subseteq (\text{col}(V^r))^\perp, \text{col}(\Delta)\subseteq \text{col}(U^r))^\perp\},\label{eq-sub2}
\end{align}
\end{subequations}
where $\text{row}(\Delta)\in \R^{d_2}$ and $\text{col}(\Delta)\in \R^{d_1}$ denote the row space and column space of the matrix $\Delta$, respectively. When the sub-matrices $(U^r,V^r)$ are clear from the context, the shorthand notation $\Aa^r$ and $\Bb^r$ are adopted instead. The definitions of $\Aa^r$ and $\Bb^r$ have been first used in \cite{agarwal2012fast,negahban2011estimation} to investigate low-rank estimation problems without measurement errors, in order to verify the decomposability of the nuclear norm, that is, for any pair of matrices $\Theta\in \Aa^r$ and $\Theta'\in \Bb^r$, it holds that $\normm{\Theta+\Theta'}_*=\normm{\Theta}_*+\normm{\Theta'}_*$, implying that the nuclear norm is decomposable with respect to the subspaces $\Aa^r$ and $\Bb^r$.



Still consider the singular decomposition $\Theta^*=UDV^\top$. For any positive number $\eta>0$ to be chosen, we define the set corresponding to $\Theta^*$:
\begin{equation}\label{eq-thresh}
K_\eta:=\{j\in\{1,2,\cdots,d\}\big||\sigma_j(\Theta^*)|>\eta\}.
\end{equation}
By the above notations, the matrix $U^{|K_\eta|}$ (resp., $V^{|K_\eta|}$) represents the $d_1\times |K_\eta|$ (resp., the $d_2\times |K_\eta|$) orthogonal matrix consisting of the singular vectors corresponding to the top $|K_\eta|$ singular values of $\Theta^*$.
Then the matrix $\Theta^*_{K_\eta^c}:=\Pi_{B^{|K_\eta|}}(\Theta^*)$ has rank at most $d-|K_\eta|$, with singular values $\{\sigma_j(\Theta^*
), j\in K_\eta^c\}$. Moreover, recall the true parameter $\Theta^*\in \B_q(R_q)$. Then the cardinality of $K_\eta$ and the approximation error $\normm{\Theta^*_{K_\eta^c}}_*$ can both be upper bounded. Indeed, it follows directly from a standard argument (see, e.g., \cite{negahban2011estimation}) that
\begin{equation}\label{eq-thresh-bound}
|K_\eta|\leq \eta^{-q}R_q\quad \mbox{and}\quad \normm{\Theta^*_{K_\eta^c}}_*\leq\eta^{1-q}R_q.
\end{equation}
Essentially, the cardinality of $K_\eta$ plays the role of the effective rank with a suitably chosen value for $\eta$ under the near low-rank assumption, as will be clarified by the proof of Theorem 1; see Remark \ref{rmk-stat}(iii).

With these notations, we now state a useful technical lemma that shows, for the true parameter matrix $\Theta^*$ and any matrix $\Theta$, we can decompose $\Delta:=\Theta-\Theta^*$ as the sum of two matrices $\Delta'$ and $\Delta''$ such that the
rank of $\Delta'$ and the nuclear norm of $\Delta''$ are both not too large. Similar results have been previously established in \cite{negahban2011estimation}, and we provide it here for the sake of completeness with the proof omitted.

\begin{Lemma}\label{lem-decom}
For a positive integer $r\in \{1, 2,\cdots,d\}$, let $U^r\in\R^{d_1\times r}$ and $V^r\in\R^{d_2\times r}$ be matrices consisting of the top $r$ left and right singular vectors in singular values of $\Theta^*$, respectively. Let $\Theta\in \R^{d_1\times d_2}$ be an arbitrary matrix. Then for the error matrix $\Delta:=\Theta-\Theta^*$, there exists a decomposition $\Delta=\Delta'+\Delta''$ such that:
\begin{enumerate}[{\rm (i)}]
\item the matrix $\Delta'$ satisfies that $\emph{rank}(\Delta')\leq 2r$;
\item moreover, suppose that $\Theta=\hat{\Theta}$ is a global optimum of the optimization problem \eqref{eq-mine}. Then if $\lambda\geq 2\phi(\Q,\sigma_\epsilon)\sqrt{\frac{\max(d_1,d_2)}{n}}$, the nuclear norm of $\Delta''$ is bounded as
\begin{equation}\label{eq-lem1-cone}
\normm{\Delta''}_*\leq 3\normm{\Delta'}_*+4\sum_{j=r+1}^d\sigma_j(\Theta^*).
\end{equation}
\end{enumerate}
\end{Lemma}

\subsection{Statistical recovery bounds}

The next theorem establishes recovery bounds for the error matrix $\hat{\Theta}-\Theta^*$ in terms of Frobenius and nuclear norms.

\begin{Theorem}\label{thm-stat}
Let $R_q>0$ and $\omega>0$ be positive numbers such that $\Theta^*\in \B_q(R_q)\cap \Su$. Let $\hat{\Theta}$ be a global optimum of the optimization problem \eqref{eq-mine}.
Suppose that the surrogate matrices $({\hat{\Gamma},\hat{\Upsilon}})$ satisfy the deviation condition \eqref{eq-devia}, and that $\hat{\Gamma}$ satisfies the \emph{RSC} condition (cf. \eqref{eq-rsc}) with
\begin{equation}\label{eq-thm1-rsc}
\tau\leq \frac{\phi(\Q,\sigma_\epsilon)}{\omega}\sqrt{\frac{\max (d_1,d_2)}{n}}.
\end{equation}
Assume that $\lambda$ is chosen to satisfy
\begin{equation}\label{eq-thm1-lambda}
\lambda\geq  2\phi(\Q,\sigma_\epsilon)\sqrt{\frac{\max (d_1,d_2)}{n}}.
\end{equation}
Then we have that
\begin{equation}\label{eq-l2bound}
\normm{\hat{\Theta}-\Theta^*}_F^2\leq 544R_q\left(\frac{\lambda}{\alpha_1}\right)^{2-q},
\end{equation}
\begin{equation}\label{eq-l1bound}
\normm{\hat{\Theta}-\Theta^*}_*\leq (4+32\sqrt{17})R_q\left(\frac{\lambda}{\alpha_1}\right)^{1-q}.
\end{equation}
\end{Theorem}

\begin{Remark}\label{rmk-stat}
{\rm (i)} Note that there are actually two parameters involved in the optimization problem \eqref{eq-mine}, namely the regularization parameter $\lambda$ and the radius of the side constraint $\omega$. On the simulated data, as the true parameter $\Theta^*$ is known beforehand, these two parameters can be chosen accordingly as we did in simulations (cf. Supplementary material). In real data analysis where $\Theta^*$ is unknown, it seems that the side constraint $\normm{\Theta^*}_*\leq \omega$ is quite restrictive. Nevertheless, Theorem \ref{thm-stat} provide some heuristic to obtain the scale for the radius $\omega$. Specifically, it follows from assumptions \eqref{eq-thm1-rsc} and \eqref{eq-thm1-lambda} that the relation $\lambda\geq 2\tau\omega$ holds. Thus one can use methods such as cross-validation to tune these two parameters based on this relation. The choice on $\lambda$ and $\omega$ in Theorem \ref{thm-algo} can also be settled in this way.

{\rm (ii)} Theorem \ref{thm-stat} shows that the recovery bound on the squared Frobenius norm for all the global solutions of the nonconvex optimization problem \eqref{eq-mine} scale as $\normm{\hat{\Theta}-\Theta^*}_F^2=O(\lambda^{2-q}R_q)$. When $\lambda$ is chosen as $\lambda=\Omega\left(\sqrt{\frac{\max(d_1,d_2)}{n}}\right)$, and given $R_q\left(\frac{\max(d_1,d_2)}{n}\right)^{1-q/2}=o(1)$, this recovery bound indicates that the estimator $\hat{\Theta}$ is statistically consistent.

{\rm (iii)}  When $q=0$, the underlying parameter $\Theta^*$ is of exact low-rank with $\emph{rank}(\Theta^*)\leq R_0$. More generally, Theorem \ref{thm-stat} provides the recovery bound when $\Theta^*\in \B_q(R_q)$ with $q\in [0,1]$. As the low-rankness of $\Theta^*$ is measured via the matrix $\ell_q$-norm and larger values means higher rank, \eqref{eq-l2bound} reveals that the rate of the recovery bound slows down along with $q$ increasing to 1. In addition, Theorem \ref{thm-stat} also offers some insight into the effective rank of a near low-rank matrix $\Theta^*$. Specifically, when $\lambda$ is chosen as $\lambda=\Omega\left(\sqrt{\frac{\max(d_1,d_2)}{n}}\right)$, and the threshold $\eta$ in \eqref{eq-thresh} is set to $\eta=\frac{\lambda}{\alpha_1}$ as we did in the proof of Theorem \ref{thm-stat}, the cardinality of the set $K_\eta$ serves as the effective rank under the near low-rank situation. This special value is chosen to provide a trade-off between the estimation and approximation errors for a near low-rank matrix. As $\lambda$ decays to 0 with the sample size $n$ increasing, this effective rank will increase (cf. \eqref{eq-thresh-bound}), a fact which implies that as more samples are obtained, it is possible to estimate more smaller singular values of a near low-rank matrix.

{\rm (v)} It is worth noting that due to the nonconvexity of the optimization problem \eqref{eq-mine},  finding a global solution is always difficult. Even worse, there may be many global optima located far apart. However, Theorem \ref{thm-stat} guarantees that under suitable conditions, the second issue is not significant any more. Specifically, suppose that there are two optima $\hat{\Theta}$ and $\hat{\Theta}'$, then it follows from Theorem \ref{thm-stat} and the triangle inequality that
\begin{equation*}
\normm{\hat{\Theta}-\hat{\Theta}'}_F\leq \normm{\hat{\Theta}-\Theta^*}_F+\normm{\hat{\Theta}'-\Theta^*}_F\leq 2c_0\lambda^{1-q/2}R_q^{1/2},
\end{equation*}
where $c_0$ is a positive constant. Hence, when $\lambda$ is chosen as $\lambda=\Omega\left(\sqrt{\frac{\max(d_1,d_2)}{n}}\right)$, and provided that $R_q\left(\frac{\max(d_1,d_2)}{n}\right)^{1-q/2}=o(1)$, all global optima must lie within an $\ell_2$-ball with radius decaying to 0. The first issue can also be handled in our setting, as we shall see in Theorem \ref{thm-algo} that the proximal gradient method converges to a local solution in polynomial time which is close enough to the set of all global optima, and thus acts as a good approximation of the global solution.
\end{Remark}

\subsection{Computational algorithms}\label{sec-comp}

Due to the nonconvexity of problem \eqref{eq-mine}, numerical methods may terminated in local optima which may be far away from any global optima. Hence, though Theorem \ref{thm-stat} provides statistical guarantee for all global optima, from the point of practice, we still need to resolve the problem on how to obtain such a global solution, or at least a near-global solution in polynomial-time. The following result tells us that under reasonable conditions, the proximal gradient method \cite{nesterov2007gradient} converges linearly to a local solution which is very near to any global solution.

We now apply the proximal gradient method \cite{nesterov2007gradient} to solve the proposed nonconvex optimization problem \eqref{eq-mine} and then establish the linear convergence result. Recall the loss function $\loss(\Theta)=\frac{1}{2}\inm{\hat{\Gamma}\Theta}{\Theta}-\inm{\hat{\Upsilon}}{\Theta}$ and the optimization objective function $\obj(\Theta)=\loss(\Theta)+\lambda\normm{\Theta}_*$. The gradient of the loss function takes the form $\nabla\loss(\Theta)=\hat{\Gamma}\Theta-\hat{\Upsilon}$. Then it is easy to see that the optimization objective function consists of a differentiable but nonconvex function and a nonsmooth but convex function (i.e., the nuclear norm). The proximal gradient method proposed in \cite{nesterov2007gradient} is applied to \eqref{eq-mine} to obtain a sequence of iterates $\{\Theta^t\}_{t=0}^\infty$ as
\begin{equation}\label{eq-algo-pga}
\Theta^{t+1}\in \argmin_{\Theta\in \Su}\left\{\frac{1}{2}\normm{\Theta-\left(\Theta^t-\frac{\nabla\loss(\Theta^t)}{v}\right)}_F^2+\frac{\lambda}{v}\normm{\Theta)}_*\right\},
\end{equation}
where $\frac1v$ is the step size.

Recall the feasible region $\Su=\{\Theta\in \R^{d_1\times d_2}\big|\normm{\Theta}_*\leq \omega\}$. According to discussions on the constrained proximal gradient update \cite[Appendix C.1]{loh2015regularized}, the next iterate $\Theta^{t+1}$ can be generated through the following three procedures.
\begin{enumerate}[(1)]
\item First solve the unconstrained optimization problem
\begin{equation*}
\hat{\Theta}^t\in \argmin_{\Theta\in \R^{d_1\times d_2}}\left\{\frac{1}{2}\normm{\Theta-\left(\Theta^t-\frac{\nabla\loss(\Theta^t)}{v}\right)}_F^2+\frac{\lambda}{v}\normm{\Theta)}_*\right\}.
\end{equation*}
\item If $\normm{\Theta^t}_*\leq \omega$, let $\Theta^{t+1}=\hat{\Theta}^t$.
\item Otherwise, if $\normm{\Theta^t}_*>\omega$, solve the constrained optimization problem
\begin{equation*}
\Theta^{t+1}\in \argmin_{\normm{\Theta}_*\leq \omega}\left\{\frac{1}{2}\normm{\Theta-\left(\Theta^t-\frac{\nabla\loss(\Theta^t)}{v}\right)}_F^2\right\}.
\end{equation*}
\end{enumerate}

Before we state our main computational result that the sequence generated by \eqref{eq-algo-pga} converges geometrically to a small neighborhood of any global solution $\hat{\Theta}$, some notations are needed first to simplify the expositions. Recall the RSC and RSM conditions in \eqref{eq-rsc} and \eqref{eq-rsm}, respectively, and that the true underlying parameter $\Theta^*\in \B_q(R_q)$ (cf. \eqref{eq-lqball}).  Let $\hat{\Theta}$ be a global solution of the optimization problem \eqref{eq-mine}. Then unless otherwise specified, we define
\begin{align}
&\bar{\epsilon}_{\text{stat}}:=8\lambda^{-\frac{q}{2}}R_q^{\frac12}\left(\sqrt{2}\normm{\hat{\Theta}-\Theta^*}_F+\lambda^{1-\frac{q}2}R_q^{\frac12}\right),\label{bar-epsilon}\\
&\kappa:= \left\{1-\frac{\alpha_1}{8v}+\frac{256\tau\lambda^{-q}R_q}{\alpha_1}\right\}\left\{1-\frac{256\tau \lambda^{-\frac{q}2}R_q}{\alpha_1}\right\}^{-1},\label{lem-bound-kappa}\\
&\xi:= \tau\left\{\frac{\alpha_1}{8v}+\frac{512\tau\lambda^{\frac{q}2}R_q}{\alpha_1}+5\right\}\left\{1-\frac{256\tau\lambda^{-\frac{q}2}R_q}{\alpha_1}\right\}^{-1}.\label{lem-bound-xi}
\end{align}

For a given number $\delta>0$ and an integer $T>0$ such that
\begin{equation}\label{lem-cone-De1}
\obj(\Theta^t)-\obj(\hat{\Theta})\leq \delta, \quad \forall\ t\geq T,
\end{equation}
define
\begin{equation*}
\epsilon(\delta):=2\min\left(\frac{\delta}{\lambda},\omega\right).
\end{equation*}

With this setup, we now state our main result on computational guarantee as follows.
\begin{Theorem}\label{thm-algo}
Let $R_q>0$ and $\omega>0$ be positive numbers such that $\Theta^*\in \B_q(R_q)\cap \Su$. Let $\hat{\Theta}$ be a global solution of the optimization problem \eqref{eq-mine}.
Suppose that $\hat{\Gamma}$ satisfies the \emph{RSC} and \emph{RSM }conditions (cf. \eqref{eq-rsc} and \eqref{eq-rsm}) with
\begin{equation}\label{eq-thm2-rsc}
\tau\leq \frac{\phi(\Q,\sigma_\epsilon)}{\omega}\sqrt{\frac{\max (d_1,d_2)}{n}}.
\end{equation}
Let $\{\Theta^t\}_{t=0}^\infty$ be a sequence of iterates generated via \eqref{eq-algo-pga} with an initial point $\Theta^0$ and $v\geq \max\{4\alpha_1, \alpha_2\}$. Assume that $\lambda$ is chosen to satisfy
\begin{equation}\label{eq-thm2-lambda}
\lambda\geq  \max\left\{\left(\frac{128\tau R_q}{\alpha_1}\right)^{1/q} , 4\phi(\Q,\sigma_\epsilon)\sqrt{\frac{\max (d_1,d_2)}{n}}\right\}.
\end{equation}
Then for any tolerance $\delta^*\geq\frac{8\xi}{1-\kappa}\bar{\epsilon}_{\emph{stat}}^2$ and any iteration $t\geq T(\delta^*)$, we have that
\begin{equation}\label{thm2-error}
\normm{\Theta^t-\hat{\Theta}}_F^2\leq \frac{4}{\alpha_1}\left(\delta^*+\frac{{\delta^*}^2}{2\tau\omega^2}+2\tau\bar{\epsilon}_{\emph{stat}}^2\right),
\end{equation}
where
\begin{equation*}
\begin{aligned}
T(\delta^*)&:=\log_2\log_2\left(\frac{\omega\lambda}{\delta^*}\right)\left(1+\frac{\log 2}{\log(1/\kappa)}\right) +\frac{\log((\obj(\Theta^0)-\obj(\hat{\Theta}))/\delta^*)}{\log(1/\kappa)},
\end{aligned}
\end{equation*}
and $\bar{\epsilon}_{\emph{stat}}$, $\kappa$, $\xi$ are defined in
\eqref{bar-epsilon}-\eqref{lem-bound-xi}, respectively.
\end{Theorem}

\begin{Remark}\label{rmk-algo}
{\rm (i)} Theorem \ref{thm-algo} establishes the upper bound on the squared Frobenius norm between $\Theta^t$ at iteration $t$, which can be easily obtained via the proposed method in polynomial time, and any global solution $\hat{\Theta}$, which is rather difficult to compute in general. Note from \eqref{thm2-error} that the optimization error $\normm{\Theta^t-\hat{\Theta}}_F$ depends on the statistical error $\normm{\hat{\Theta}-\Theta^*}_F$ with constants ignored, and thus Theorem \ref{thm-algo} guarantees that the proximal gradient method produces a local solution lying in a small neighborhood of all the global optima, and essentially performs as well as any global solution of the nonconvex problem \eqref{eq-mine}, in the sense of statistical error. Therefore, Theorem \ref{thm-algo} provides the insight that for nonconvex estimation problems, though the global solution cannot be obtained generally, a near-global solution can still be found and serves as a good candidate under certain regularity conditions.

{\rm (ii)} Theorem \ref{thm-algo} also provides the upper bound on the number of iteration counts required in a logarithmic scale to solve the nonconvex problem \eqref{eq-mine}. It should be noted that the established geometric convergence is not ensured to an arbitrary precision, but only to an accuracy associated with statistical precision. As has been pointed in \cite{agarwal2012fast} that, this is a sensible implication, since in statistical settings, there is no point to solving optimization problems beyond the statistical precision. Thus the near-global solution which the proximal gradient method converged to is the best we can hope from the point of statistical computing.

{\rm (iii)} In the case when $q=0$, the parameter $\Theta^*$ is of exact low-rank with $\emph{rank}(\Theta^*)\leq R_0$. More generally, Theorem \ref{thm-algo} establishes the linear convergence rate when $\Theta^*\in \B_q(R_q)$ with  $q\in [0,1]$ and indicates some important differences between the case of exact low-rank and that of near low-rank, a finding which has also been discussed in \cite{agarwal2012fast,li2020sparse} for sparse linear regression. Specifically, for the exact low-rank case (i.e., $q=0$), the parameter $\bar{\epsilon}_\emph{stat}$ defined in \eqref{bar-epsilon} actually does not contain the second term $\lambda^{1-q}R_q$ (To see this, turn to the proof of Lemma D.1 (cf. inequality (D.8)) and note the fact that $\normm{\Pi_{\Bb^r}(\Theta^*)}_*=0$. We did not point out this case explicitly in the theorem for the sake of compactness of the article.)  Thus the optimization error $\normm{\Theta^t-\hat{\Theta}}_F$ only relies on the statistical error $\normm{\hat{\Theta}-\Theta^*}_F$ in this exact low-rank case. In contrast, apart from the squared statistical error $\normm{\hat{\Theta}-\Theta^*}_F^2$, the squared optimization error \eqref{thm2-error} in the case when $q\in (0,1]$ also has an additional term $\lambda^{2-2q}R_q^2$ (to see this, take the square of \eqref{bar-epsilon} with some calculation), which appears due to the statistical nonidentifiability over the matrix $\ell_q$-ball, and it is no larger than $\normm{\hat{\Theta}-\Theta^*}_F^2$ with high probability.
\end{Remark}

\color{black}\section{Consequences for errors-in-variables regression}\label{sec-conse}

Both Theorems \ref{thm-stat} and \ref{thm-algo} are intrinsically
deterministic as mentioned before. Consequences for specific models requires probabilistic discussions to verify that the core conditions are satisfied, namely, the RSC/RSM conditions (cf. \eqref{eq-rsc}/\eqref{eq-rsm})
and the deviation condition (cf. \eqref{eq-devia}). Verification of these conditions involves nontrivial analysis on concentration inequalities and random matrix theory, and serves as one of the main contributions of this article.

We now turn to establish consequences for different cases of additive noise and missing data by applying Theorems \ref{thm-stat} and \ref{thm-algo}. The following propositions are needed first to facilitate the establishment.

\subsection{Additive noise case}

In the additive noise case, set $\Sigma_z=\Sigma_x+\Sigma_w$, $\sigma_z^2=\normm{\Sigma_x}_\text{op}^2+\sigma_w^2$ for notational simplicity. Then define
\begin{subequations}
\begin{align}
\tau_{\text{add}}&=\lambda_{\min}(\Sigma_x)\max\left(\frac{d_2^2(\normm{\Sigma_x}_\text{op}^2+\sigma_w^2)^2}{\lambda_{\min}^2(\Sigma_x)},\frac{d_2(\normm{\Sigma_x}_\text{op}^2+\sigma_w^2)}{\lambda_{\min}(\Sigma_x)}\right)\frac{2\max (d_1,d_2)+\log(\min(d_1,d_2))}{n},\label{tau-add}\\
\phi_{\text{add}}&=\sqrt{\lambda_{\max}(\Sigma_z)}(\sigma_\epsilon+\omega\sigma_w).\label{phi-add}
\end{align}
\end{subequations}

\begin{Proposition}[RSC/RSM conditions, additive noise case]\label{prop-add-rs}
Let $\tau_{\emph{add}}$ be given in \eqref{tau-add}. In the additive noise case, there exist universal positive constants $(c_0, c_1)$ such that the matrix $\hat{\Gamma}_{\emph{add}}$ satisfies the \emph{RSC} and \emph{RSM }conditions (cf. \eqref{eq-rsc} and \eqref{eq-rsm}) with parameters $\alpha_1=\frac{\lambda_{\min}(\Sigma_x)}{2}$, $\alpha_2=\frac{3\lambda_{\max}(\Sigma_x)}{2}$, and $\tau=c_0\tau_{\emph{add}}$, with probability at least\\
$1-2\exp\left(-c_1n\min\left(\frac{\lambda_{\min}^2(\Sigma_x)}{d_2^2(\normm{\Sigma_x}_\emph{op}^2+\sigma_w^2)^2},\frac{\lambda_{\min}(\Sigma_x)}{d_2(\normm{\Sigma_x}_\emph{op}^2+\sigma_w^2)}\right)+\log d_2\right)$.
\end{Proposition}

\begin{Proposition}[Deviation condition, additive noise case]\label{prop-add-devia}
Let $\phi_{\emph{add}}$ be given in \eqref{phi-add}. In the additive noise case, there exist universal positive constants $(c_0,c_1,c_2)$ such that the deviation condition (cf. \eqref{eq-devia}) holds with parameter $\phi(\Q,\sigma_\epsilon)= c_0\phi_{\emph{add}}$, with probability at least $1-c_1\exp(-c_2\log (\max(d_1,d_2))$.
\end{Proposition}

\begin{Remark}
The prefactor $\phi_\emph{add}$ in the deviation bound actually serves as the integrate noise term in the additive noise case. Obviously, $\phi_\emph{add}$ scales up with $\sigma_\epsilon$ and $\sigma_w$—two sources of errors
in this setting. Then when $\sigma_\epsilon=0$, corresponding to the case that there is no observation error $\epsilon$, it holds that $\phi_\emph{add}\neq 0$ because of additive noise in the observed covariates; while when $\sigma_w=0$, corresponding to that the covariates are correctly obtained, it still holds that $\phi_\emph{add}\neq 0$ due to the observation error introduced by $\epsilon$, and the result then reduces to the known finding for clean covariates in low-rank estimation problems \cite{negahban2011estimation}.
\end{Remark}

Now we are ready to state statistical and computational consequences for the multi-response regression with additive noise. The conclusions follow by applying Propositions \ref{prop-add-rs} and \ref{prop-add-devia} on Theorems \ref{thm-stat} and \ref{thm-algo}, respectively, and so the proofs are omitted.

\begin{Corollary}\label{corol-add-stat}
Let $R_q>0$ and $\omega>0$ be positive numbers such that $\Theta^*\in \B_q(R_q)\cap \Su$. Let $\hat{\Theta}$ be a global optimum of the optimization problem \eqref{eq-mine} with $(\hat{\Gamma}_{\emph{add}},\hat{\Upsilon}_{\emph{add}})$ given by \eqref{sur-add} in place of $(\hat{\Gamma},\hat{\Upsilon})$. Then there exist universal positive constants $c_i\ (i=0,1,2,3,4)$ such that
if $\tau_{\emph{add}}\leq c_0\frac{\phi_\emph{add}}{\omega}\sqrt{\frac{\max (d_1,d_2)}{n}}$ and $\lambda$ is chosen to satisfy
$\lambda\geq  c_1\phi_{\emph{add}}\sqrt{\frac{\max (d_1,d_2)}{n}}$, then it holds with probability at least $1-2\exp\left(-c_2n\min\left(\frac{\lambda_{\min}^2(\Sigma_x)}{d_2^2(\normm{\Sigma_x}_\emph{op}^2+\sigma_w^2)^2},\frac{\lambda_{\min}(\Sigma_x)}{d_2(\normm{\Sigma_x}_\emph{op}^2+\sigma_w^2)}\right)+\log d_2\right)-c_3\exp(-c_4\log (\max(d_1,d_2))$
that
\begin{equation*}
\normm{\hat{\Theta}-\Theta^*}_F^2\leq 544R_q\left(\frac{2\lambda}{\lambda_{\min}(\Sigma_x)}\right)^{2-q},
\end{equation*}
\begin{equation*}
\normm{\hat{\Theta}-\Theta^*}_*\leq (4+32\sqrt{17})R_q\left(\frac{2\lambda}{\lambda_{\min}(\Sigma_x)}\right)^{1-q}.
\end{equation*}
\end{Corollary}
\begin{Remark}\label{rmk-corol1}
In the special case when $\Sigma_x=\sigma_x^2\I_{d_1}$, the integrate noise term takes the value $\phi_{\emph{add}}=\sqrt{\sigma_x^2+\sigma_w^2}(\sigma_\epsilon+\omega\sigma_w)$. Then if $\lambda$ is chosen as $\lambda=\Omega\left(\phi_{\text{add}}\sqrt{\frac{\max(d_1,d_2)}{n}}\right)$, the squared Frobenius error bound takes the form $\normm{\hat{\Theta}-\Theta^*}_F^2=O\left(\left(\frac{\sqrt{\sigma_x^2+\sigma_w^2}(\sigma_\epsilon+\omega\sigma_w)}{\sigma_x^2}\right)^{2-q}\left(\frac{\max(d_1,d_2)}{n}\right)^{1-q/2}R_q\right)$. When $\Theta^*$ is an exact low-rank matrix, i.e., $q=0$ and $\emph{rank}(\Theta^*)\leq R_0$, this error bound further reduces to $\normm{\hat{\Theta}-\Theta^*}_F^2=O\left(\left(\frac{\sqrt{\sigma_x^2+\sigma_w^2}(\sigma_\epsilon+\omega\sigma_w)}{\sigma_x^2}\right)^2\frac{\max(d_1,d_2)R_0}{n}\right)$. The scaling can be explained as that the squared error is proportional to the noise $\left(\frac{\sqrt{\sigma_x^2+\sigma_w^2}(\sigma_\epsilon+\omega\sigma_w)}{\sigma_x^2}\right)^2$ and the quantity $\max(d_1,d_2)R_0$ which measures the degrees of freedom of a matrix with size $d_1\times d_2$ and rank $R_0$ up to constant factors. Note that if there exists no constraint on $\Theta^*$, then since a $d_1\times d_2$ matrix possesses $d_1d_2$ free parameters, the statistical convergence rate will at best take the order $\normm{\hat{\Theta}-\Theta^*}_F^2=O\left(\left(\frac{\sqrt{\sigma_x^2+\sigma_w^2}(\sigma_\epsilon+\omega\sigma_w)}{\sigma_x^2}\right)^2\frac{d_1d_2}{n}\right)$. Nevertheless, when $\Theta^*$ is supposed to be exact low-rank with $R_0\ll \min(d_1,d_2)$, then the nuclear norm regularized estimator gains a substantially faster statistical convergence rate. This benefit of the nuclear norm regularization in reducing the freedom has also been discussed in \cite{negahban2011estimation}. In addition, we also note that it requires the quantity $\max(d_1,d_2)$ to be large enough in high-dimensional problems to ensure that the error probability to be small. In this case, the conclusions (cf. Corollaries \ref{corol-add-stat}--\ref{corol-mis-algo}) can be held with overwhelming probability, and such a high-dimensional scaling is just the main focus of this work.
\end{Remark}

Define
\begin{align*}
&\kappa_\text{add}:= \left\{1-\frac{\lambda_{\min}(\Sigma_x)}{16v}+\frac{512\tau_{\text{add}}\lambda^{-q}R_q}{\lambda_{\min}(\Sigma_x)}\right\}\left\{1-\frac{512\tau_{\text{add}}\lambda^{-\frac{q}2}R_q}{\lambda_{\min}(\Sigma_x)}\right\}^{-1},
\\
&\xi_\text{add}:= \tau_{\text{add}}\left\{\frac{\lambda_{\min}(\Sigma_x)}{16v}+\frac{1024\tau_{\text{add}}\lambda^{\frac{q}2}R_q}{\lambda_{\min}(\Sigma_x)}+5\right\}\left\{1-\frac{512\tau_{\text{add}}\lambda^{-\frac{q}2}R_q}{\lambda_{\min}(\Sigma_x)}\right\}^{-1}.
\end{align*}

\begin{Corollary}\label{corol-add-algo}
Let $R_q>0$ and $\omega>0$ be positive numbers such that $\Theta^*\in \B_q(R_q)\cap \Su$. Let $\hat{\Theta}$ be a global solution of the optimization problem \eqref{eq-mine} with $(\hat{\Gamma}_{\emph{add}},\hat{\Upsilon}_{\emph{add}})$ given by \eqref{sur-add} in place of $(\hat{\Gamma},\hat{\Upsilon})$. Let $\{\Theta^t_{\emph{add}}\}_{t=0}^\infty$ be a sequence of iterates generated via \eqref{eq-algo-pga} with an initial point $\Theta^0_{\emph{add}}$ and $v\geq \max(2\lambda_{\min}(\Sigma_x), \frac{3}{2}\lambda_{\max}(\Sigma_x))$.
Then there exist universal positive constants $c_i\ (i=0,1,2,3,4,5)$ such that, if
$\tau_{\emph{add}}\leq c_0\frac{\phi_{\emph{add}}}{\omega}\sqrt{\frac{\max (d_1,d_2)}{n}}$ and $\lambda$ is chosen to satisfy
\begin{equation*}
\lambda\geq  \max\left\{\left(\frac{c_1\tau_{\emph{add}}R_q}{\lambda_{\min}(\Sigma_x)}\right)^{1/q}, c_2\phi_{\emph{add}}\sqrt{\frac{\max (d_1,d_2)}{n}}\right\},
\end{equation*}
then for any tolerance $\delta^*\geq\frac{8\xi_{\emph{add}}}{1-\kappa_{\emph{add}}}\bar{\epsilon}_{\emph{stat}}^2$ and any iteration $t\geq T(\delta^*)$, it holds with probability at least
\begin{flalign*}
&1-2\exp\left(-c_3n\min\left(\frac{\lambda_{\min}^2(\Sigma_x)}{d_2^2(\normm{\Sigma_x}_\emph{op}^2+\sigma_w^2)^2},\frac{\lambda_{\min}(\Sigma_x)}{d_2(\normm{\Sigma_x}_\emph{op}^2+\sigma_w^2)}\right)+\log d_2\right)\qquad \\
&-c_4\exp(-c_5\log (\max(d_1,d_2))
\end{flalign*}
that
\begin{equation*}
\normm{\Theta^t_{\emph{add}}-\hat{\Theta}}_F^2\leq \frac{8}{\lambda_{\min}(\Sigma_x)}\left(\delta^*+\frac{{\delta^*}^2}{2\tau_{\emph{add}}\omega^2}+2\tau_{\emph{add}}\bar{\epsilon}_{\emph{stat}}^2\right),
\end{equation*}
where
\begin{equation*}
\begin{aligned}
T(\delta^*)&:=\log_2\log_2\left(\frac{\omega\lambda}{\delta^*}\right)\left(1+\frac{\log 2}{\log(1/\kappa_{\emph{add}})}\right) +\frac{\log((\obj(\Theta^0_{\emph{add}})-\obj(\hat{\Theta}))/\delta^*)}{\log(1/\kappa_{\emph{add}})},
\end{aligned}
\end{equation*}
and $\bar{\epsilon}_{\emph{stat}}$ is given in
\eqref{bar-epsilon}.
\end{Corollary}

\subsection{Missing data case}

In the missing data case, define a matrix $M\in \R^{d_1\times d_1}$ satisfying $M_{ij}=(1-\rho)^2$ for $i\not =j$ and $M_{ij}=1-\rho$ for $i=j$. Let $\otimes$ and $\oslash$ denote element-wise multiplication and division, respectively, and set $\Sigma_z=\Sigma_x\otimes M$. Then define
\begin{subequations}
\begin{align}
\tau_{\text{mis}}&=\lambda_{\min}(\Sigma_x)\max\left(\frac{1}{(1-\rho)^4}\frac{d_2^2\normm{\Sigma_x}_\text{op}^4}{\lambda_{\min}^2(\Sigma_x)},\frac{1}{(1-\rho)^2}\frac{d_2\normm{\Sigma_x}_\text{op}^2}{\lambda_{\min}(\Sigma_x)}\right)\frac{2\max(d_1,d_2)+\log(\min(d_1,d_2))}{n},\label{tau-mis}\\
\phi_{\text{mis}}&=\frac{\sqrt{\lambda_{\max}(\Sigma_z)}}{1-\rho}\left(\frac{\omega}{1-\rho}\normm{\Sigma_x}_\text{op}+\sigma_\epsilon\right).\label{phi-mis}
\end{align}
\end{subequations}
As with $\phi_\text{add}$, the prefactor $\phi_\text{mis}$ in the deviation bound also serves as the integrate noise term in the missing data case. 

\begin{Proposition}[RSC/RSM conditions, missing data case]\label{prop-mis-rs}
Let $\tau_{\emph{mis}}$ be given in \eqref{tau-mis}. In the missing data case, there exist universal positive constants $(c_0, c_1)$ such that the matrix $\hat{\Gamma}_{\emph{mis}}$ satisfies the \emph{RSC} and \emph{RSM} conditions (cf. \eqref{eq-rsc} and \eqref{eq-rsm}) with parameters $\alpha_1=\frac{\lambda_{\min}(\Sigma_x)}{2}$, $\alpha_2=\frac{3\lambda_{\max}(\Sigma_x)}{2}$, and $\tau=c_0\tau_{\emph{mis}}$, with probability at least\\
$1-2\exp\left(-c_1n\min\left((1-\rho)^4\frac{\lambda_{\min}^2(\Sigma_x)}{d_2^2\normm{\Sigma_x}_\emph{op}^4},(1-\rho)^2\frac{\lambda_{\min}(\Sigma_x)}{d_2\normm{\Sigma_x}_\emph{op}^2}\right)+\log d_2\right)$.
\end{Proposition}

\begin{Proposition}[Deviation condition, missing data case]\label{prop-mis-devia}
Let $\phi_{\emph{mis}}$ be given in \eqref{phi-mis}. In the missing data case, there exist universal positive constants $(c_0,c_1,c_2)$ such that the deviation condition (cf. \eqref{eq-devia}) holds with parameter $\phi(\Q,\sigma_\epsilon)=c_0\phi_{\emph{mis}}$, with probability at least $1-c_1\exp(-c_2\max(d_1,d_2))$.
\end{Proposition}


Now we are at the stage to state concrete statistical and computational properties for the multi-response regression with missing data. The conclusions follow by applying Propositions \ref{prop-mis-rs} and \ref{prop-mis-devia} on Theorems \ref{thm-stat} and \ref{thm-algo}, respectively, and so the proofs are omitted.

\begin{Corollary}\label{corol-mis-stat}
Let $R_q>0$ and $\omega>0$ be positive numbers such that $\Theta^*\in \B_q(R_q)\cap \Su$. Let $\hat{\Theta}$ be a global optimum of the optimization problem \eqref{eq-mine} with $(\hat{\Gamma}_{\emph{mis}},\hat{\Upsilon}_{\emph{mis}})$ given by \eqref{sur-mis} in place of $(\hat{\Gamma},\hat{\Upsilon})$. Then there exist universal positive constants $c_i\ (i=0,1,2,3,4)$ such that
if $\tau_{\emph{mis}}\leq c_0\frac{\phi_\emph{mis}}{\omega}\sqrt{\frac{\max (d_1,d_2)}{n}}$, and $\lambda$ is chosen to satisfy
$\lambda\geq  c_1\phi_{\emph{mis}}\sqrt{\frac{\max (d_1,d_2)}{n}}$, then it holds with probability at least $1-2\exp\left(-c_2n\min\left((1-\rho)^4\frac{\lambda_{\min}^2(\Sigma_x)}{d_2^2\normm{\Sigma_x}_\emph{op}^4},(1-\rho)^2\frac{\lambda_{\min}(\Sigma_x)}{d_2\normm{\Sigma_x}_\emph{op}^2}\right)+\log d_2\right)-c_3\exp(-c_4\max( d_1,d_2))$
that
\begin{equation*}
\normm{\hat{\Theta}-\Theta^*}_F^2\leq 544R_q\left(\frac{2\lambda}{\lambda_{\min}(\Sigma_x)}\right)^{2-q},
\end{equation*}
\begin{equation*}
\normm{\hat{\Theta}-\Theta^*}_*\leq (4+32\sqrt{17})R_q\left(\frac{2\lambda}{\lambda_{\min}(\Sigma_x)}\right)^{1-q}.
\end{equation*}
\end{Corollary}


Define
\begin{align*}
&\kappa_\text{mis}:= \left\{1-\frac{\lambda_{\min}(\Sigma_x)}{16v}+\frac{512\tau_{\text{mis}}R_q\lambda^{-q}}{\lambda_{\min}(\Sigma_x)}\right\}\left\{1-\frac{512\tau_{\text{mis}}R_q\lambda^{-\frac{q}2}}{\lambda_{\min}(\Sigma_x)}\right\}^{-1},
\\
&\xi_\text{mis}:= \tau_{\text{mis}}\left\{\frac{\lambda_{\min}(\Sigma_x)}{16v}+\frac{1024\tau_{\text{mis}}R_q\lambda^{\frac{q}2}}{\lambda_{\min}(\Sigma_x)}+5\right\}\left\{1-\frac{512\tau_{\text{mis}}R_q\lambda^{-\frac{q}2}}{\lambda_{\min}(\Sigma_x)}\right\}^{-1}.
\end{align*}

\begin{Corollary}\label{corol-mis-algo}
Let $R_q>0$ and $\omega>0$ be positive numbers such that $\Theta^*\in \B_q(R_q)\cap \Su$. Let $\hat{\Theta}$ be a global solution of the optimization problem \eqref{eq-mine} with $(\hat{\Gamma}_{\emph{mis}},\hat{\Upsilon}_{\emph{mis}})$ given by \eqref{sur-mis} in place of $(\hat{\Gamma},\hat{\Upsilon})$. Let $\{\Theta^t_{\emph{mis}}\}_{t=0}^\infty$ be a sequence of iterates generated via \eqref{eq-algo-pga} with an initial point $\Theta^0_{\emph{mis}}$ and $v\geq \max(2\lambda_{\min}(\Sigma_x), \frac{3}{2}\lambda_{\max}(\Sigma_x))$.
Then there exist universal positive constants $c_i\ (i=0,1,2,3,4,5)$ such that, if
$\tau_{\emph{mis}}\leq c_0\frac{\phi_{\emph{mis}}}{\omega}\sqrt{\frac{\max (d_1,d_2)}{n}}$ and $\lambda$ is chosen to satisfy
\begin{equation*}
\lambda\geq  \max\left\{\left(\frac{c_1\tau_{\emph{mis}}R_q}{\lambda_{\min}(\Sigma_x)}\right)^{1/q}, c_2\phi_{\emph{mis}}\sqrt{\frac{\max (d_1,d_2)}{n}}\right\},
\end{equation*}
then for any tolerance $\delta^*\geq\frac{8\xi_{\emph{mis}}}{1-\kappa_{\emph{mis}}}\bar{\epsilon}_{\emph{stat}}^2$ and any iteration $t\geq T(\delta^*)$, it holds with probability at least
\begin{flalign*}
&1-2\exp\left(-c_3n\min\left((1-\rho)^4\frac{\lambda_{\min}^2(\Sigma_x)}{d_2^2\normm{\Sigma_x}_\emph{op}^4},(1-\rho)^2\frac{\lambda_{\min}(\Sigma_x)}{d_2\normm{\Sigma_x}_\emph{op}^2}\right)+\log d_2\right)\qquad \\
&-c_4\exp(-c_5\max( d_1,d_2))
\end{flalign*}
that
\begin{equation*}
\normm{\Theta^t_{\emph{mis}}-\hat{\Theta}}_F^2\leq \frac{8}{\lambda_{\min}(\Sigma_x)}\left(\delta^*+\frac{{\delta^*}^2}{2\tau_{\emph{mis}}\omega^2}+2\tau_{\emph{mis}}\bar{\epsilon}_{\emph{stat}}^2\right),
\end{equation*}
where
\begin{equation*}
\begin{aligned}
T(\delta^*)&:=\log_2\log_2\left(\frac{\omega\lambda}{\delta^*}\right)\left(1+\frac{\log 2}{\log(1/\kappa_{\emph{mis}})}\right) +\frac{\log((\obj(\Theta^0_{\emph{mis}})-\obj(\hat{\Theta}))/\delta^*)}{\log(1/\kappa_{\emph{mis}})},
\end{aligned}
\end{equation*}
and $\bar{\epsilon}_{\emph{stat}}$ is given in
\eqref{bar-epsilon}.
\end{Corollary}

\section{Simulations}\label{sec-simul}

In this section, we implement several numerical experiments on the multi-response errors-in-variables regression model to illustrate our main theoretical results. The following simulations will be performed with the loss function $\loss_n$ corresponding to the additive noise and missing data cases, respectively, and with the the nuclear norm regularizer. All numerical experiments are performed in MATLAB R2013a and executed on a personal desktop (Intel Core i7-6700, 2.80 GHz, 16.00 GB of RAM).

The simulated data are generated as follows. Specifically, the true parameter is generated as a square matrix $\Theta^*\in \R^{d\times d}$, and we consider the exact low-rank case as an instance with $\text{rank}(\Theta^*)=r=10$. Explicitly, let $\Theta^*=AB^\top$, where $A, B\in\R^{d\times r}$ consist of i.i.d. $\N(0,1)$ entries. Then we generate i.i.d. true covariates $X_{i\cdot}\sim \N(0,\I_d)$ and the noise term $e\sim \N(0,(0.1)^2\I_n)$. The data $y$ are then generated according to \eqref{eq-model-matrix}. The corrupted term is set to $W_{i\cdot}\sim \N(0,(0.2)^2\I_d)$ and $\rho=0.2$ for the additive noise and missing data cases, respectively. The problem sizes $d$ and $n$ will be specified according to specific experiments. The data are generated randomly for 100 times.

For all simulations, the parameter in \eqref{eq-mine} is set as $\lambda=\sqrt{\frac{d}{n}}$, and $\omega=1.1\normm{\Theta^*}_*$ to ensure the feasibility of $\Theta^*$. Iteration \eqref{eq-algo-pga} is then carried out with $v=2\lambda_{\max}(\Sigma_x)$ and the initial point $\Theta^0$ being a zero matrix. Performance of the estimator $\hat{\Theta}$ is characterized by the relative error $\normm{\hat{\Theta}-\Theta^*}_F/\normm{\Theta^*}_F$ and is illustrated by averaging across the 100 numerical results.

\begin{figure}[htbp]
\centering
\subfigure[]{
\includegraphics[width=0.48\textwidth]{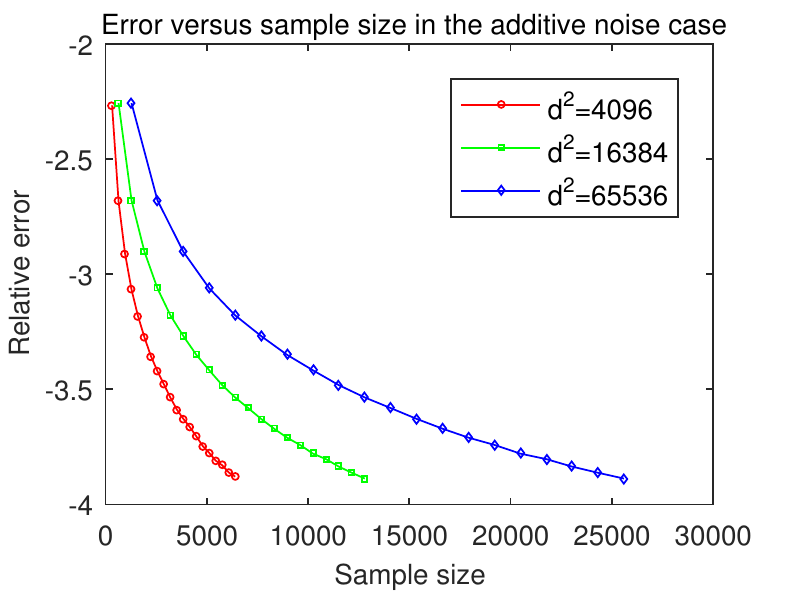}}
\subfigure[]{
\includegraphics[width=0.48\textwidth]{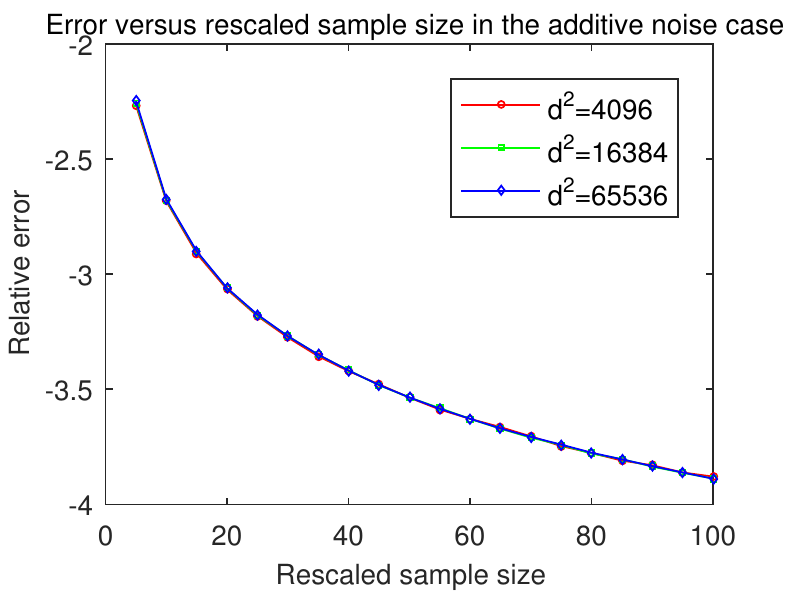}}
\caption{Statistical consistency for multi-response regression with additive error.}
\label{f-stat-add}
\end{figure}

\begin{figure}[htbp]
\centering
\subfigure[]{
\includegraphics[width=0.48\textwidth]{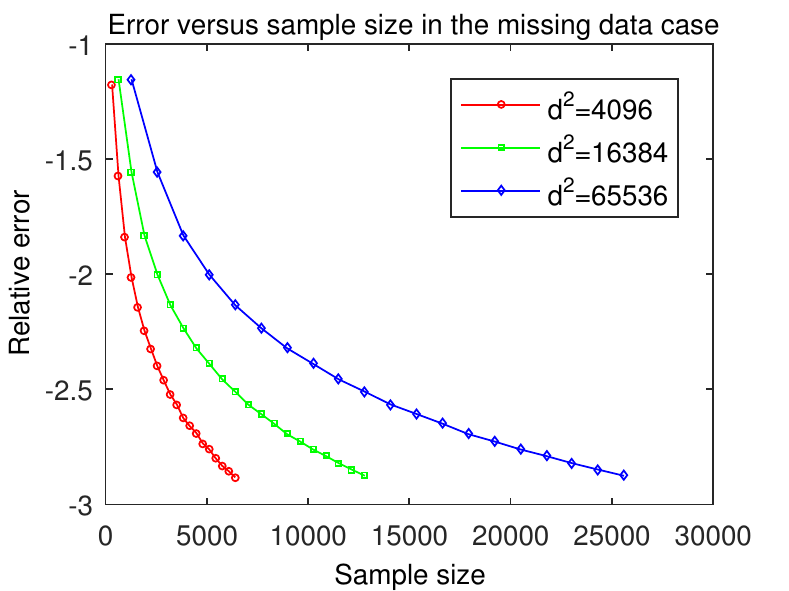}}
\subfigure[]{
\includegraphics[width=0.48\textwidth]{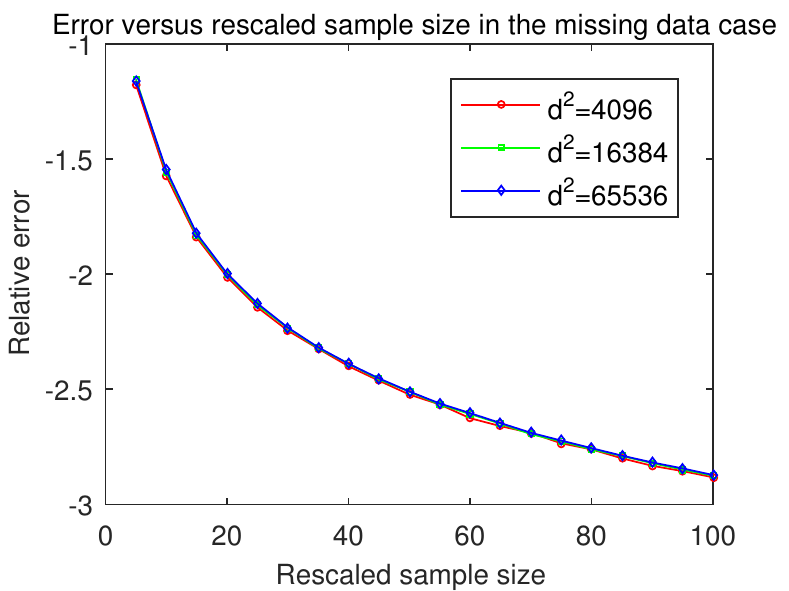}}
\caption{Statistical consistency for multi-response regression with missing data.}
\label{f-stat-mis}
\end{figure}

The first experiment is performed to demonstrate the statistical guarantee for multi-response regression in additive noise and missing data cases, respectively. Fig. \ref{f-stat-add}(a) plots the relative error on a logarithmic scale versus the sample size $n$ for three different matrix dimensions $d\in \{64,128,256\}$ in the additive noise case. For each matrix dimension, as the sample size increases, the relative error decays to zero, implying the statistical consistency. However, more samples are needed to estimate a larger matrix, as is illustrated by the rightward shift of the curves along with the dimension $d$ increasing. Fig. \ref{f-stat-add}(b) shows the same set of numerical results as in Fig. \ref{f-stat-add}(a), but now the relative error is plotted versus the rescaled sample size $n/d$. We can see from Fig. \ref{f-stat-add}(b) that the three curves perfectly  match with one another under different matrix dimensions in the sense that these curves are now unrecognizable, a fact which has been predicted by Corollary \ref{corol-add-stat}. Hence, Fig. \ref{f-stat-add} shows that the quantity $n/d$ actually plays the role of the effective sample size in this high-dimensional setting. Similar results on the statistical consistency and the ``stacking'' phenomenon under the effective sample size for the missing data case are displayed in Fig. \ref{f-stat-mis}.

\begin{figure}[htbp]
\centering
\subfigure[]{
\includegraphics[width=0.48\textwidth]{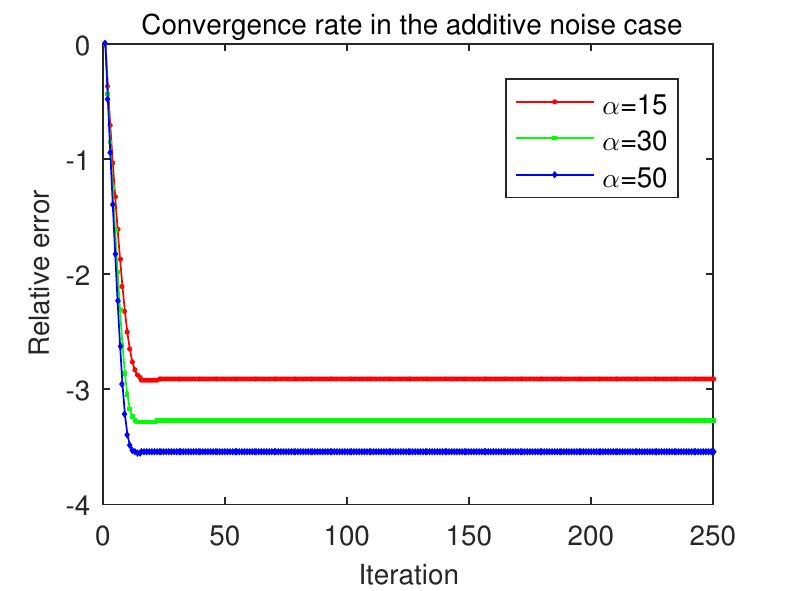}}
\subfigure[]{
\includegraphics[width=0.48\textwidth]{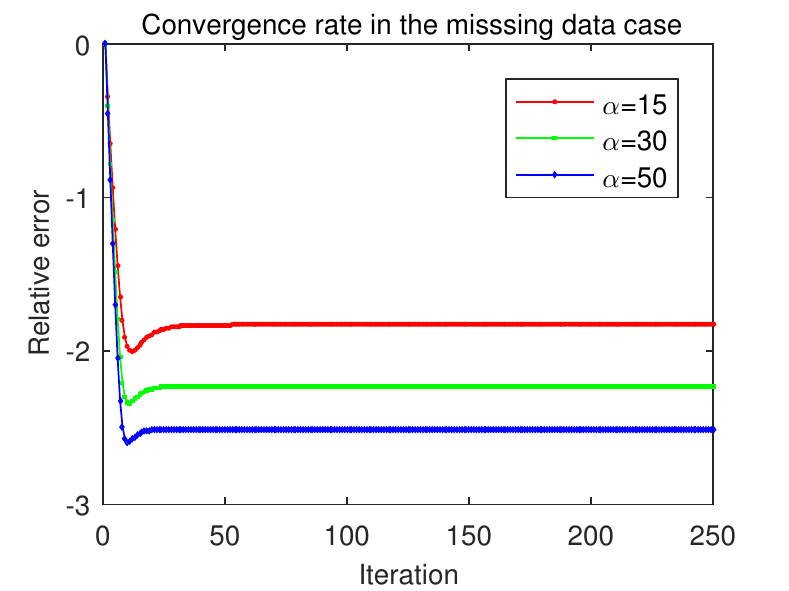}}
\caption{Algorithmic convergence rate for multi-response errors-in-variables regression.}
\label{f-algo}
\end{figure}

The second experiment is designed to illustrate the algorithmic linear convergence rate in additive noise and missing data cases, respectively. We have explored the performance for a broad range of dimensions $d$ and $n$, and the results are fairly consistent across these choices. Thus we here only report results for $d=128$ and a range of the sample sizes $n=\lceil\alpha d\rceil$ with $\alpha\in\{15,30,50\}$. In the additive noise case, Fig. \ref{f-algo}(a) shows that for the three sample sizes, the algorithm reveals exact linear convergence rate. As the sample size becomes larger, the convergence speed turns faster and achieves a more accurate estimation level. Fig. \ref{f-algo}(b) depicts analogous results to Fig. \ref{f-algo}(a) in the missing data case.

\section{Conclusion}\label{sec-concl}

In this work, we investigated the low-rank estimation problem in high-dimensional multi-response errors-in-variables regression. A nonconvex error-corrected estimator was constructed. Then under the (near) low-rank situation, the theoretical statistical and computational properties for global solutions of the nonconvex estimator were analysed. In the statistical aspect, we established recovery bound for all the global optima to imply statistical consistency. In the computational aspect, the proximal gradient method was applied to solve the optimization problem and was proved to converge to a near-global solution in polynomial time. Probabilistic consequences were obtained for specific errors-in-variables models via verifying the regularity conditions. Theoretical consequences were illustrated by several numerical experiments on multi-response errors-in-variables regression models.

Future directions of research include modifying the proposed method to be adaptive to the low-rankness of the underlying parameter, to overcome the obstacle that the hyper-parameters $q$ and the $\ell_q$-radius are usually unknown aforehand in implementation. Strategies such as combining or aggregating will  be a powerful tool; see e.g., \cite{wang2014adaptive}. In addition, though we assumed that the error covariance $\Sigma_w$ or the missing probability $\rho$ is known, it is still an open question on how to construct suitable estimators for them. Finally, it would be interesting to extend our method to study more general types of multiplicative noise or the case where data that are not Gaussian. 
\vskip 10pt

\noindent \textbf{Acknowledgements} Xin Li's work was supported in part by the National Natural Science Foundation of China (Grant No. 12201496) and the Natural Science Foundation of Shaanxi Province of China (Grant No. 2022JQ-045). Dongya Wu's work was supported in part by the National Natural Science Foundation of China (Grant No. 62103329).

\appendix

\section{Proof of Theorem \ref{thm-stat}}

\setcounter{equation}{0}
\renewcommand{\theLemma}{A.\arabic{Lemma}}
\renewcommand{\theequation}{A.\arabic{equation}}

Set $\hat{\Delta}:=\hat{\Theta}-\Theta^*$. By the feasibility of $\Theta^*$ and optimality of $\hat{\Theta}$, one has that $\obj(\hat{\Theta})\leq \obj(\Theta^*)$. Then it follows from elementary algebra and the triangle inequality that
\begin{equation*}
\frac{1}{2}\inm{\hat{\Gamma}\hat{\Delta}}{\hat{\Delta}}
\leq \inm{\hat{\Upsilon}-\hat{\Gamma}\Theta^*}{\hat{\Delta}}+\lambda\normm{\Theta^*}_*-\lambda\normm{\Theta^*+\hat{\Delta}}_* \leq \inm{\hat{\Upsilon}-\hat{\Gamma}\Theta^*}{\hat{\Delta}}+\lambda\normm{\hat{\Delta}}_*.
\end{equation*}
Applying H{\"o}lder's inequality and by the deviation condition \eqref{eq-devia}, one has that
\begin{equation*}
\inm{\hat{\Upsilon}-\hat{\Gamma}\Theta^*}{\hat{\Delta}}\leq \phi(\Q,\sigma_\epsilon)\sqrt{\frac{\max(d_1,d_2)}{n}}\normm{\hat{\Delta}}_*.
\end{equation*}
Combining the above two inequalities, and noting \eqref{eq-thm1-lambda},
we obtain that
\begin{equation}\label{eq-thm12}
\inm{\hat{\Gamma}\hat{\Delta}}{\hat{\Delta}}
\leq 3\lambda\normm{\hat{\Delta}}_*.
\end{equation}
Applying the RSC condition \eqref{eq-rsc} to the left-hand side of \eqref{eq-thm12}, yields that
\begin{equation}\label{eq-thm16}
\alpha_1\normm{\hat{\Delta}}_F^2-\tau\normm{\hat{\Delta}}_*^2\leq 3\lambda\normm{\hat{\Delta}}_*.
\end{equation}
On the other hand, by assumptions \eqref{eq-thm1-rsc} and \eqref{eq-thm1-lambda} and noting the fact that $\normm{\hat{\Delta}}_*\leq \normm{\Theta^*}_*+\normm{\hat{\Theta}}_*\leq 2\omega$,
the left-hand side of \eqref{eq-thm16} is lower bounded as
\begin{equation*}
\alpha_1\normm{\hat{\Delta}}_F^2-\tau\normm{\hat{\Delta}}_*^2
\geq
\alpha_1\normm{\hat{\Delta}}_F^2-2\tau\omega\normm{\hat{\Delta}}_*
\geq
\alpha_1\normm{\hat{\Delta}}_F^2-\lambda\normm{\hat{\Delta}}_*.
\end{equation*}
Combining this inequality with \eqref{eq-thm16}, one has that
\begin{equation}\label{eq-thm17}
\alpha_1\normm{\hat{\Delta}}_F^2
\leq 4\lambda\normm{\hat{\Delta}}_*.
\end{equation}
Then it follows from Lemma \ref{lem-decom} that there exists a matrix $\hat{\Delta}'$ such that
\begin{equation}\label{eq-thm13}
\normm{\hat{\Delta}}_*\leq 4\normm{\hat{\Delta}'}_*+4\sum_{j=r+1}^d\sigma_j(\Theta^*)\leq 4\sqrt{2r}\normm{\hat{\Delta}'}_F+4\sum_{j=r+1}^d\sigma_j(\Theta^*),
\end{equation}
where $\text{rank}(\hat{\Delta}')\leq 2r$ with $r$ to be chosen later, and the second inequality is due to the fact that $\normm{\hat{\Delta}'}_*\leq \sqrt{2r}\normm{\hat{\Delta}'}_F$.
Combining \eqref{eq-thm17} and \eqref{eq-thm13}, we obtain that
\begin{equation*}
\alpha_1\normm{\hat{\Delta}}_F^2
\leq 16\lambda\left(\sqrt{2r}\normm{\hat{\Delta}'}_F+\sum_{j=r+1}^d\sigma_j(\Theta^*)\right)
\leq 16\lambda\left(\sqrt{2r}\normm{\hat{\Delta}}_F+\sum_{j=r+1}^d\sigma_j(\Theta^*)\right).
\end{equation*}
Then it follows that
\begin{equation}\label{eq-thm1-rbound}
\normm{\hat{\Delta}}_F^2\leq \frac{512r\lambda^2+32\alpha_1\lambda\sum_{j=r+1}^d\sigma_j(\Theta^*)}{\alpha_1^2}.
\end{equation}
Recall the set $K_\eta$ defined in \eqref{eq-thresh} and set $r=|K_\eta|$. Combining \eqref{eq-thm1-rbound} with \eqref{eq-thresh-bound} and setting $\eta=\frac{\lambda}{\alpha_1}$, we arrive at \eqref{eq-l2bound}. Moreover, it follows from \eqref{eq-thm13} that \eqref{eq-l1bound} holds. The proof is complete.

\section{Proof of Theorem \ref{thm-algo}}

\setcounter{equation}{0}
\setcounter{Lemma}{0}
\renewcommand{\theLemma}{B.\arabic{Lemma}}
\renewcommand{\theequation}{B.\arabic{equation}}

Before providing the proof of Theorem \ref{thm-algo}, we need several useful lemmas first.

\begin{Lemma}\label{lem-cone}
Suppose that the conditions of Theorem \ref{thm-algo} are satisfied, and
that there exists a pair $(\delta, T)$ such that \eqref{lem-cone-De1} holds.
Then for any iteration $t\geq T$, it holds that
\begin{equation*}
\begin{aligned}
\normm{\Theta^t-\hat{\Theta}}_*&\leq 4\sqrt{2}\lambda^{-\frac{q}{2}}R_q^{\frac12}\normm{\Theta^t-\hat{\Theta}}_F+\bar{\epsilon}_{\emph{stat}}+\epsilon(\delta).
\end{aligned}
\end{equation*}
\end{Lemma}
\begin{proof}
We first prove that if $\lambda\geq 4\phi(\Q,\sigma_\epsilon)\sqrt{\frac{\max(d_1,d_2)}{n}}$, then for any $\Theta\in \Su$ satisfying
\begin{equation}\label{lem-cone-De2}
\obj(\Theta)-\obj(\Theta^*)\leq \delta,
\end{equation}
it holds that
\begin{equation}\label{lem-cone-1}
\begin{aligned}
\normm{\Theta-\Theta^*}_*&\leq 4\sqrt{2}\lambda^{-\frac{q}{2}}R_q^{\frac12}\normm{\Theta-\Theta^*}_F+4\lambda^{1-q}R_q+2\min\left(\frac{\delta}{\lambda},\omega\right).
\end{aligned}
\end{equation}
Set $\Delta:=\Theta-\Theta^*$. From \eqref{lem-cone-De2}, one has that
\begin{equation*}
\loss(\Theta^*+\Delta)+\lambda\normm{\Theta^*+\Delta}_*\leq \loss(\Theta^*)+\lambda\normm{\Theta^*}_*+\delta.
\end{equation*}
Then subtracting $\inm{\nabla\loss(\Theta^*)}{\Delta}$ from both sides of the former inequality and recalling the formulation of $\loss(\cdot)$, we obtain that
\begin{equation}\label{lem-cone-2}
\frac{1}{2}\inm{\hat{\Gamma}\Delta}{\Delta}+\lambda\normm{\Theta^*+\Delta}_*-\lambda\normm{\Theta^*}_*\leq -\inm{\hat{\Gamma}\Theta^*-\hat{\Upsilon}}{\Delta}+\delta.
\end{equation}
We now claim that
\begin{equation}\label{lem-cone-claim}
\lambda\normm{\Theta^*+\Delta}_*-\lambda\normm{\Theta^*}_*\leq \frac{\lambda}{2}\normm{\Delta}_*+\delta.
\end{equation}
In fact, combining \eqref{lem-cone-2} with the RSC condition \eqref{eq-rsc} and H{\"o}lder's inequality, one has that
\begin{equation*}
\frac{1}{2}\{\alpha_1\normm{\Delta}_F^2-\tau\normm{\Delta}_*^2\}+\lambda\normm{\Theta^*+\Delta}_*-\lambda\normm{\Theta^*}_*\leq \normm{\hat{\Upsilon}-\hat{\Gamma}\Theta^*}_{\text{op}}\normm{\Delta}_*+\delta.
\end{equation*}
This inequality, together with the deviation condition \eqref{eq-devia} and the assumption that $\lambda\geq 4\phi(\Q,\sigma_\epsilon)\sqrt{\frac{\max(d_1,d_2)}{n}}$, implies that
\begin{equation*}
\frac{1}{2}\{\alpha_1\normm{\Delta}_F^2-\tau\normm{\Delta}_*^2\}+\lambda\normm{\Theta^*+\Delta}_*-\lambda\normm{\Theta^*}_*\leq \frac{\lambda}{4}\normm{\Delta}_*+\delta.
\end{equation*}
Noting the facts that $\alpha_1>0$ and that $\normm{\Delta}_*\leq \normm{\Theta^*}_*+\normm{\Theta}_*\leq 2\omega$, one arrives at \eqref{lem-cone-claim} by combining assumptions \eqref{eq-thm2-rsc} and \eqref{eq-thm2-lambda}. On the other hand, it follows from Lemma \ref{lem-decom}(i) that there exists two matrices $\Delta'$ and $\Delta''$ such that $\Delta=\Delta'+\Delta''$, where $\text{rank}(\Delta')\leq 2r$ with $r$ to be chosen later. Recall the definitions of $\Aa^r$ and $\Bb^r$ given respectively in \eqref{eq-sub1} and \eqref{eq-sub2}. Then the decomposition $\Theta^*=\Pi_{\Aa^r}(\Theta^*)+\Pi_{\Bb^r}(\Theta^*)$ holds. This equality, together with the triangle inequality and Lemma \ref{lem-decom}(i) as well as \eqref{eq-lem12}, implies that
\begin{equation*}
\begin{aligned}
\normm{\Theta}_*&=\normm{(\Pi_{\Aa^r}(\Theta^*)+\Delta'')+(\Pi_{\Bb^r}(\Theta^*)+\Delta')}_*\\
&\geq \normm{\Pi_{\Aa^r}(\Theta^*)+\Delta'')}_*-\normm{\Pi_{\Bb^r}(\Theta^*)+\Delta')}_*\\
&\geq
\normm{\Pi_{\Aa^r}(\Theta^*)}_*+\normm{\Delta''}_*-\{\normm{\Pi_{\Bb^r}(\Theta^*)}_*+\normm{\Delta'}_*\}.
\end{aligned}
\end{equation*}
Consequently, we have
\begin{equation}\label{lem-cone-8}
\begin{aligned}
\normm{{\Theta}^*}_*-\normm{\Theta}_*
&\leq
\normm{{\Theta}^*}_*-\normm{\Pi_{\Aa^r}(\Theta^*)}-\normm{\Delta''}_*+\{\normm{\Pi_{\Bb^r}(\Theta^*)}_*+\normm{\Delta'}_*\}\\
&\leq
2\normm{\Pi_{\Bb^r}(\Theta^*)}_*+\normm{\Delta'}_*-\normm{\Delta''}_*.
\end{aligned}
\end{equation}
Combining \eqref{lem-cone-8} and \eqref{lem-cone-claim} and noting the fact that $\normm{\Pi_{\Bb^r}(\Theta^*)}_*=\sum_{j=r+1}^{d}\sigma_j(\Theta^*)$, one has that $0\leq \frac{3\lambda}{2}\normm{\Delta'}_*-\frac{\lambda}{2}\normm{\Delta''}_*+2\lambda \sum_{j=r+1}^{d}\sigma_j(\Theta^*)+\delta$, and consequently, $\normm{\Delta''}_*\leq 3\normm{\Delta'}_*+4\sum_{j=r+1}^{d}\sigma_j(\Theta^*)+\frac{2\delta}{\lambda}$. Using the trivial bound $\normm{\Delta}_*\leq
2\omega$, one has that
\begin{equation}\label{lem-cone-9}
\normm{\Delta}_*\leq 4\sqrt{2r}\normm{\Delta}_F+4\sum_{j=r+1}^d\sigma_j(\Theta^*)+2\min\left(\frac{\delta}{\lambda},\omega\right).
\end{equation}
Recall the set $K_\eta$ defined in \eqref{eq-thresh} and set $r=|K_\eta|$. Combining \eqref{lem-cone-9} with \eqref{eq-thresh-bound} and setting $\eta=\lambda$, we arrive at \eqref{lem-cone-1}. We now verify that \eqref{lem-cone-De2} is held by the matrix $\hat{\Theta}$ and $\Theta^t$, respectively. Since $\hat{\Theta}$ is the optimal solution, it holds that $\obj(\hat{\Theta})\leq \obj(\Theta^*)$, and by assumption \eqref{lem-cone-De1}, it holds that $\obj(\Theta^t)\leq \obj(\hat{\Theta})+\delta\leq \obj(\Theta^*)+\delta$. Consequently, it follows from \eqref{lem-cone-1} that
\begin{equation*}
\begin{aligned}
\normm{\hat{\Theta}-\Theta^*}_*&\leq 4\sqrt{2}\lambda^{-\frac{q}{2}}R_q^{\frac12}\normm{\hat{\Theta}-\Theta^*}_F+4\lambda^{1-q}R_q,\\
\normm{\Theta^t-\Theta^*}_*&\leq 4\sqrt{2}\lambda^{-\frac{q}{2}}R_q^{\frac12}\normm{\Theta^t-\Theta^*}_F+4\lambda^{1-q}R_q+2\min\left(\frac{\delta}{\lambda},\omega\right).
\end{aligned}
\end{equation*}
By the triangle inequality, we then arrive at that
\begin{equation*}
\begin{aligned}
\normm{\Theta^t-\hat{\Theta}}_*
&\leq \normm{\hat{\Theta}-\Theta^*}_*+\normm{\Theta^t-\Theta^*}_*\\
&\leq 4\sqrt{2}\lambda^{-\frac{q}{2}}R_q^{\frac12}(\normm{\hat{\Theta}-\Theta^*}_F+\normm{\Theta^t-\Theta^*}_F)+8\lambda^{1-q}R_q+2\min\left(\frac{\delta}{\lambda},\omega\right)\\
&\leq 4\sqrt{2}\lambda^{-\frac{q}{2}}R_q^{\frac12}\normm{\Theta^t-\hat{\Theta}}_F+\bar{\epsilon}_{\text{stat}}+\epsilon(\delta).
\end{aligned}
\end{equation*}
The proof is complete.
\end{proof}

\begin{Lemma}\label{lem-Tphi-bound}
Suppose that the conditions of Theorem \ref{thm-algo} are satisfied and that there exists a pair $(\Delta, T)$ such that \eqref{lem-cone-De1} holds.
Then for any iteration $t\geq T$, we have that
\begin{align}
&\loss(\hat{\Theta})-\loss(\Theta^t)-\inm{\nabla\loss(\Theta^t)}{\hat{\Theta}-\Theta^t}\geq -\tau(\bar{\epsilon}_{\emph{stat}}+\epsilon(\delta))^2,\label{lem-bound-T1}\\
&\obj(\Theta^t)-\obj(\hat{\Theta})\geq \frac{\alpha_1}{4}\normm{\Theta^t-\hat{\Theta}}_F^2-\tau(\bar{\epsilon}_{\emph{stat}}+\epsilon(\delta))^2, \label{lem-bound-T2}\\
&\obj(\Theta^t)-\obj(\hat{\Theta})\leq \kappa^{t-T}(\obj(\Theta^T)-\obj(\hat{\Theta}))+\frac{2\xi}{1-\kappa}(\bar{\epsilon}_{\emph{stat}}^2+\epsilon^2(\delta)).\label{lem-bound-0}
\end{align}
\end{Lemma}
\begin{proof}
By the RSC condition \eqref{eq-rsc}, one has that
\begin{equation}\label{lem-bound-1}
\loss(\Theta^t)-\loss(\hat{\Theta})-\inm{\nabla\loss(\hat{\Theta})}{\Theta^t-\hat{\Theta}}\geq \frac12\left\{\alpha_1\normm{\Theta^t-\hat{\Theta}}_F^2-\tau\normm{\Theta^t-\hat{\Theta}}_*^2\right\}.
\end{equation}
It then follows from Lemma \ref{lem-cone} and the assumption that
$\lambda\geq \left(\frac{128\tau R_q}{\alpha_1}\right)^{1/q}$ that
\begin{equation*}
\loss(\hat{\Theta})-\loss(\Theta^t)-\inm{\nabla\loss(\Theta^t)}{\hat{\Theta}-\Theta^t}
\geq \frac12\left\{\alpha_1\normm{\Theta^t-\hat{\Theta}}_F^2-\tau\normm{\Theta^t-\hat{\Theta}}_*^2\right\}
\geq -\tau(\bar{\epsilon}_{\text{stat}}+\epsilon(\delta))^2,
\end{equation*}
which establishes \eqref{lem-bound-T1}. Furthermore, by the convexity of $\normm{\cdot}_*$, one has that
\begin{equation}\label{lem-bound-2}
\lambda\normm{\Theta^t}_*-\lambda\normm{\hat{\Theta}}_*-\inm{\nabla\left\{\lambda\normm{\hat{\Theta}}_*\right\}}{\Theta^t-\hat{\Theta}} \geq 0,
\end{equation}
and by the first-order optimality condition for $\hat{\Theta}$, one has that
\begin{equation}\label{lem-bound-3}
\inm{\nabla\obj(\hat{\Theta})}{\Theta^t-\hat{\Theta}}\geq 0.
\end{equation}
Combining \eqref{lem-bound-1}, \eqref{lem-bound-2} and \eqref{lem-bound-3}, we obtain that
\begin{equation*}
\obj(\Theta^t)-\obj(\hat{\Theta})\geq \frac12\left\{\alpha_1\normm{\Theta^t-\hat{\Theta}}_F^2-\tau\normm{\Theta^t-\hat{\Theta}}_*^2\right\}.
\end{equation*}
Then applying Lemma \ref{lem-cone} to bound the term $\normm{\Theta^t-\hat{\Theta}}_*^2$ and noting the assumption that $\lambda\geq \left(\frac{128\tau R_q}{\alpha_1}\right)^{1/q}$,
we arrive at \eqref{lem-bound-T2}.
Now we turn to prove \eqref{lem-bound-0}. Define
\begin{equation*}
\obj_t(\Theta):=\loss(\Theta^t)+\inm{\nabla{\loss(\Theta^t)}}{\Theta-\Theta^t}+\frac{v}{2}\normm{\Theta-\Theta^t}_F^2+\lambda\normm{\Theta}_*,
\end{equation*}
which is the objective function minimized over the feasible region $\Su=\{\Theta\big|\normm{\Theta}_*\leq \omega\}$ at iteration count $t$. For any $a\in [0,1]$, it is easy to check that the matrix $\Theta_a=a\hat{\Theta}+(1-a)\Theta^t$ belongs to $\Su$  due to the convexity of $\Su$. Since $\Theta^{t+1}$ is the optimal solution of the optimization problem \eqref{eq-algo-pga}, we have that
\begin{equation*}
\begin{aligned}
\obj_t(\Theta^{t+1}) &\leq \obj_t(\Theta_a)=\loss(\Theta^t)+\inm{\nabla{\loss(\Theta^t)}}{\Theta_a-\Theta^t} +\frac{v}{2}\|\Theta_a-\Theta^t\|_2^2+\lambda\normm{\Theta_a}_*\\
&\leq \loss(\Theta^t)+ \inm{\nabla\loss(\Theta^t)}{a\hat{\Theta}-a\Theta^t} +\frac{va^2}{2}\normm{\hat{\Theta}-\Theta^t}_F^2+a\lambda\normm{\hat{\Theta}}_*+(1-a)\lambda\normm{\Theta^t}_*,
\end{aligned}
\end{equation*}
where the last inequality is from the convexity of $\normm{\cdot}_*$.
Then by \eqref{lem-bound-T1}, one sees that
\begin{equation}\label{lem-bound-4}
\begin{aligned}
\obj_t(\Theta^{t+1})
&\leq (1-a)\loss(\Theta^t)+a\loss(\hat{\Theta})+a\tau(\bar{\epsilon}_{\text{stat}}+\epsilon(\delta))^2\\
&\quad +\frac{va^2}{2}\normm{\hat{\Theta}-\Theta^t}_F^2+a\lambda\normm{\hat{\Theta}}_*+(1-a)\lambda\normm{\Theta^t}_*\\
&\leq \obj(\Theta^t)-a(\obj(\Theta^t)-\obj(\hat{\Theta}))+\tau(\bar{\epsilon}_{\text{stat}}+\epsilon(\delta))^2+\frac{va^2}{2}\normm{\hat{\Theta}-\Theta^t}_F^2.
\end{aligned}
\end{equation}
Applying the RSM condition \eqref{eq-rsm} on the matrix $\Theta^{t+1}-\Theta^t$ with some algebra, we have by assumption $v\geq \alpha_2$ that
\begin{equation*}
\begin{aligned}
\loss(\Theta^{t+1})-\loss(\Theta^t)-\inm{\nabla\loss(\Theta^t)}{\Theta^{t+1}-\Theta^t}
&\leq \frac12\left\{\alpha_2\normm{\Theta^{t+1}-\Theta^t}_F^2+\tau\normm{\Theta^{t+1}-\Theta^t}_*^2\right\}\\
&\leq \frac{v}{2}\normm{\Theta^{t+1}-\Theta^t}_F^2+\frac{\tau}{2}\normm{\Theta^{t+1}-\Theta^t}_*^2.
\end{aligned}
\end{equation*}
Adding $\lambda\normm{\Theta^{t+1}}_*$ to both sides of the former inequality, we obtain that
\begin{equation*}
\begin{aligned}
\obj(\Theta^{t+1})
&\leq \loss(\Theta^t)+\inm{\nabla\loss(\Theta^t)}{\Theta^{t+1}-\Theta^t} +\lambda\normm{\Theta^{t+1}}_*+\frac{v}{2}\normm{\Theta^{t+1}-\Theta^t}_F^2+\frac{\tau}{2}\normm{\Theta^{t+1}-\Theta^t}_*^2\\
&=
\obj_t(\Theta^{t+1})+\frac{\tau}{2}\normm{\Theta^{t+1}-\Theta^t}_*^2.
\end{aligned}
\end{equation*}
This, together with \eqref{lem-bound-4}, implies that
\begin{equation}\label{lem-bound-5}
\obj(\Theta^{t+1})\leq \obj(\Theta^t)-a(\obj(\Theta^t)-\obj(\hat{\Theta}))+\frac{va^2}{2}\normm{\hat{\Theta}-\Theta^t}_F^2+\frac{\tau}{2}\normm{\Theta^{t+1}-\Theta^t}_*^2+\tau(\bar{\epsilon}_{\text{stat}}+\epsilon(\delta))^2.
\end{equation}
Define $\Delta^t:=\Theta^t-\hat{\Theta}$. Then it follows directly that
$\normm{\Theta^{t+1}-\Theta^t}_*^2\leq (\normm{\Delta^{t+1}}_*+\normm{\Delta^t}_*)^2\leq 2\normm{\Delta^{t+1}}_*^2+2\normm{\Delta^t}_*^2$. Combining this inequality with \eqref{lem-bound-5}, one has that
\begin{equation*}
\obj(\Theta^{t+1})\leq \obj(\Theta^t)-a(\obj(\Theta^t)-\obj(\hat{\Theta}))+\frac{va^2}{2}\normm{\hat{\Theta}-\Theta^t}_F^2+\tau(\normm{\Delta^{t+1}}_*^2+\normm{\Delta^t}_*^2)+\tau(\bar{\epsilon}_{\text{stat}}+\epsilon(\delta))^2.
\end{equation*}
To simplify the notations, define $\psi:=\tau(\bar{\epsilon}_{\text{stat}}+\epsilon(\Delta))^2$, $\zeta:=\tau\lambda^{-q}R_q$ and $\delta_t:=\obj(\Theta^t)-\obj(\hat{\Theta})$. Using Lemma \ref{lem-cone} to bound the term $\normm{\Delta^{t+1}}_*^2$ and $\normm{\Delta^t}_*^2$, we arrive at that
\begin{equation}\label{lem-bound-6}
\begin{aligned}
\obj(\Theta^{t+1})
&\leq \obj(\Theta^t)-a(\obj(\Theta^t)-\obj(\hat{\Theta}))+\frac{va^2}{2}\normm{\Delta^t}_F^2+64\tau\lambda^{-q}R_q(\normm{\Delta^{t+1}}_F^2+\normm{\Delta^t}_F^2)+5\psi\\
&= \obj(\Theta^t)-a(\obj(\Theta^t)-\obj(\hat{\Theta}))+\left(\frac{va^2}{2}+64\zeta\right)\normm{\Delta^t}_F^2+64\zeta\normm{\Delta^{t+1}}_F^2+5\psi.
\end{aligned}
\end{equation}
Subtracting $\obj(\hat{\Theta})$ from both sides of \eqref{lem-bound-6}, one has by \eqref{lem-bound-T2} that
\begin{equation*}
\begin{aligned}
\delta_{t+1}
&\leq (1-a)\delta_t+\frac{2va^2+256\zeta}{\alpha_1}(\delta_t+\psi) +\frac{256\zeta}{\alpha_1}(\delta_{t+1}+\psi)+5\psi.
\end{aligned}
\end{equation*}
Setting $a=\frac{\alpha_1}{4v}\in (0,1)$, it follows from the former inequality that
\begin{equation*}
\begin{aligned}
\left(1-\frac{256\zeta}{\alpha_1}\right)\delta_{t+1}&\leq \left(1-\frac{\alpha_1}{8v}+\frac{256\zeta}{\alpha_1}\right)\delta_t +\left(\frac{\alpha_1}{8v}+\frac{512\zeta}{\alpha_1}+5\right)\psi,
\end{aligned}
\end{equation*}
or equivalently, $\delta_{t+1}\leq \kappa\delta_t+\xi(\bar{\epsilon}_{\text{stat}}+\epsilon(\delta))^2$, where $\kappa$ and $\xi$ were previously defined in \eqref{lem-bound-kappa} and \eqref{lem-bound-xi}, respectively. Finally, we conclude that
\begin{equation*}
\begin{aligned}
\Delta_t &\leq \kappa^{t-T}\Delta_T+\xi(\bar{\epsilon}_{\text{stat}}+\epsilon(\delta))^2(1+\kappa+\kappa^2+\cdots+\kappa^{t-T-1})\\
&\leq \kappa^{t-T}\Delta_T+\frac{\xi}{1-\kappa}(\bar{\epsilon}_{\text{stat}}+\epsilon(\delta))^2\leq
\kappa^{t-T}\Delta_T+\frac{2\xi}{1-\kappa}(\bar{\epsilon}_{\text{stat}}^2+\epsilon^2(\delta)).
\end{aligned}
\end{equation*}
The proof is complete.
\end{proof}

By virtue of the above lemmas, we are now ready to prove Theorem \ref{thm-algo}. The proof mainly follows the arguments in \cite{agarwal2012fast,li2020sparse}.

We first prove the following inequality:
\begin{equation}\label{thm2-1}
\obj(\Theta^t)-\obj(\hat{\Theta})\leq \delta^*,\quad \forall t\geq T(\delta^*).
\end{equation}
Divide iterations $t=0,1,\cdots$ into a sequence of disjoint epochs $[T_k,T_{k+1}]$ and define the associated sequence of tolerances $\delta_0>\delta_1>\cdots$ such that
\begin{equation*}
\obj(\Theta^t)-\obj(\hat{\Theta})\leq \delta_k,\quad \forall t\geq T_k,
\end{equation*}
as well as the corresponding error term $\epsilon_k:=2\min \left\{\frac{\delta_k}{\lambda},\omega\right\}$. The values of $\{(\delta_k,T_k)\}_{k\geq 1}$ will be decided later.
Then at the first iteration,
Lemma \ref{lem-Tphi-bound} (cf. \eqref{lem-bound-0}) is applied with $\epsilon_0=2\omega$ and $T_0=0$
to conclude that
\begin{equation}\label{thm2-T0}
\obj(\Theta^t)-\obj(\hat{\Theta})\leq \kappa^t(\obj(\Theta^0)-\obj(\hat{\Theta}))+\frac{2\xi}{1-\kappa}(\bar{\epsilon}_{\text{stat}}^2 +4\omega^2),\quad \forall t\geq T_0.
\end{equation}
Set $\delta_1:=\frac{4\xi}{1-\kappa}(\bar{\epsilon}_{\text{stat}}^2 +4\omega^2)$. Noting that $\kappa\in (0,1)$ by assumption, it follows from \eqref{thm2-T0} that for $T_1:=\lceil \frac{\log (2\delta_0/\delta_1)}{\log (1/\kappa)} \rceil$,
\begin{equation*}
\begin{aligned}
\obj(\Theta^t)-\obj(\hat{\Theta})&\leq \frac{\delta_1}{2}+\frac{2\xi}{1-\kappa}\left(\bar{\epsilon}_{\text{stat}}^2 +4\omega^2\right)=\delta_1\leq \frac{8\xi}{1-\kappa}\max\left\{\bar{\epsilon}_{\text{stat}}^2,4\omega^2\right\},\quad \forall t\geq T_1.
\end{aligned}
\end{equation*}
For $k\geq 1$, define
\begin{equation}\label{thm2-DeltaT}
\delta_{k+1}:= \frac{4\xi}{1-\kappa}(\bar{\epsilon}_{\text{stat}}^2+\epsilon_k^2)
\quad \mbox{and}\quad
T_{k+1}:= \left\lceil \frac{\log (2\delta_k/\delta_{k+1})}{\log (1/\kappa)}+T_k \right\rceil.
\end{equation}
Then Lemma \ref{lem-Tphi-bound} (cf. \eqref{lem-bound-0}) is applied to concluding that for all $t\geq T_k$,
\begin{equation*}
\obj(\Theta^t)-\obj(\hat{\Theta})\leq \kappa^{t-T_k}(\obj(\Theta^{T_k})-\obj(\hat{\Theta}))+\frac{2\xi}{1-\kappa}(\bar{\epsilon}_{\text{stat}}^2+\epsilon_k^2),
\end{equation*}
which further implies that
\begin{equation*}
\obj(\Theta^t)-\obj(\hat{\Theta})\leq \delta_{k+1}\leq \frac{8\xi}{1-\kappa}\max\{\bar{\epsilon}_{\text{stat}}^2,\epsilon_k^2\},\quad \forall t\geq T_{k+1}.
\end{equation*}
From \eqref{thm2-DeltaT}, one obtains the recursion for $\{(\delta_k,T_k)\}_{k=0}^\infty$ as follows
\begin{subequations}
\begin{align}
\delta_{k+1}&\leq \frac{8\xi}{1-\kappa}\max\{\epsilon_k^2,\bar{\epsilon}_{\text{stat}}^2\},\label{thm2-recur-Delta}\\
T_k&\leq k+\frac{\log (2^k\delta_0/\delta_k)}{\log (1/\kappa)}.\label{thm2-recur-T}
\end{align}
\end{subequations}
Then by \cite[Section 7.2]{agarwal2012supplementaryMF}, it is easy to see that \eqref{thm2-recur-Delta} implies that
\begin{equation}\label{thm2-recur2}
\delta_{k+1}\leq \frac{\delta_k}{4^{2^{k+1}}}\quad \mbox{and}\quad \frac{\delta_{k+1}}{\lambda}\leq \frac{\omega}{4^{2^k}},\quad \forall k\geq 1.
\end{equation}
Now let us show how to determine the smallest $k$ such that $\delta_k\leq \delta^*$ by using \eqref{thm2-recur2}. If we are at the first epoch, \eqref{thm2-1} is clearly held due to \eqref{thm2-recur-Delta}.
Otherwise, from \eqref{thm2-recur-T}, one sees that $\delta_k\leq \delta^*$ is held after at most $$k(\delta^*)\geq \frac{\log(\log(\omega\lambda/\delta^*)/\log 4)}{\log(2)}+1=\log_2\log_2(\omega\lambda/\delta^*)$$ epoches.
Combining the above bound on $k(\delta^*)$ with \eqref{thm2-recur-T}, one obtains
that $\obj(\Theta^t)-\obj(\hat{\Theta})\leq \delta^*$ holds for all iterations
\begin{equation*}
t\geq
\log_2\log_2\left(\frac{\omega\lambda}{\delta^*}\right)\left(1+\frac{\log 2}{\log(1/\kappa)}\right)+\frac{\log(\delta_0/\delta^*)}{\log(1/\kappa)},
\end{equation*}
which establishes \eqref{thm2-1}.
Finally, as \eqref{thm2-1} is proved, one has by \eqref{lem-bound-T2} in Lemma \ref{lem-Tphi-bound}
and the assumption that $\lambda\geq \left(\frac{128\tau R_q}{\alpha_1}\right)^{1/q}$ that, for any $t\geq T(\delta^*)$,
\begin{equation*}
\frac{\alpha_1}{4}\normm{\Theta^t-\hat{\Theta}}_F^2
\leq \obj(\Theta^t)-\obj(\hat{\Theta}) +\tau\left(\epsilon(\delta^*)+\bar{\epsilon}_{\text{stat}}\right)^2
\leq \delta^*+\tau\left(\frac{2\delta^*}{\lambda}+\bar{\epsilon}_{\text{stat}}\right)^2.
\end{equation*}
Consequently,
it follows from \eqref{eq-thm2-rsc} and \eqref{eq-thm2-lambda} that, for any $t\geq T(\delta^*)$,
\begin{equation*}
\normm{\Theta^t-\hat{\Theta}}_F^2\leq \frac{4}{\alpha_1}\left(\delta^*+\frac{{\delta^*}^2}{2\tau\omega^2}+2\tau\bar{\epsilon}_{\text{stat}}^2\right).
\end{equation*}
The proof is complete.

\color{black}\section{Proofs of Section \ref{sec-conse}}

\setcounter{equation}{0}
\setcounter{Lemma}{0}
\renewcommand{\theLemma}{C.\arabic{Lemma}}
\renewcommand{\theequation}{C.\arabic{equation}}

In this section, several important technical lemmas are provided first, which are used to verify the RSC/RSM conditions and deviation conditions for specific errors-in-variables models (cf. Propositions \ref{prop-add-rs}--\ref{prop-mis-devia}).
Some notations are needed to ease the expositions. For a symbol $x\in \{0,*,F\}$ and a positive real number $r\in \R^+$, define $\M_x(r):=\{A\in R^{d_1\times d_2}|\normm{A}_x\leq r\}$, where $\normm{A}_0$ denotes the rank of matrix A. Then define the sparse set $\K(r):=\M_0(r)\cap \M_F(1)$ and the
cone set $\C(r):=\{A\in \R^{d_1\times d_2}\big|\normm{A}_*\leq \sqrt{r}\normm{A}_F\}$. The following lemma tells us that the intersection of the matrix $\ell_1$-ball with the matrix $\ell_2$-ball can be bounded by virtue of a simpler set.

\begin{Lemma}\label{lem-inter}
For any constant $r\geq 1$, it holds that
\begin{equation*}
\M_*(\sqrt{r})\cap \M_F(1)\subseteq 2\emph{cl}\{\emph{conv}\{\M_0(r)\cap \M_F(1)\}\},
\end{equation*}
where $\emph{cl}\{\cdot\}$ and $\emph{conv}\{\cdot\}$
denote the topological closure and convex hull, respectively.
\end{Lemma}
\begin{proof}
Note that when $r> \min\{d_1,d_2\}$, this containment is trivial, since the right-hand set is equal to $\M_F(2)$ and the left-hand set is contained in $\M_F(1)$. Thus, we will assume $1\leq r\leq \min\{d_1,d_2\}$.

Let $A\in \M_*(\sqrt{r})\cap \M_F(1)$. Then it follows that $\normm{A}_*\leq \sqrt{r}$ and $\normm{A}_F\leq 1$. Consider a singular value decomposition of $A$: $A=UDV^\top$,
where $U\in \R^{d_1\times d_1}$ and $V\in \R^{d_2\times d_2}$ are orthogonal matrices, and $D\in \R^{d_1\times d_2}$ consists of $\sigma_1(D),\sigma_2(D),\cdots,\sigma_k(D)$ on the ``diagonal'' and 0 elsewhere with $k=\text{rank}(A)$. Write $D=\text{diag}(\sigma_1(D),\sigma_2(D),\cdots,\sigma_k(D),0\cdots,0)$, and use $\text{vec}(D)$ to denote the vectorized form of the matrix $D$. Then it follows that $\|\text{vec}(D)\|_1\leq \sqrt{r}$ and $\|\text{vec}(D)\|_2\leq 1$. Partition the support of $\text{vec}(D)$ into disjoint subsets $T_1,T_2,\cdots$, such that $T_1$ is the index set corresponding to the first $r$
largest elements in absolute value of $\text{vec}(D)$, $T_2$ indexes the next $r$ largest elements, and so on. Write $D_i=\text{diag}(\text{vec}(D)_{T_i})$, and $A_i=UD_iV^\top$. Then one has that $\normm{A_i}_0=\text{rank}(A_i)\leq r$ and $\normm{A_i}_F\leq 1$. Write $B_i=2A_i/\normm{A_i}_F$ and $t_i=\normm{A_i}_F/2$. Then it holds that $B_i\in 2\{\M_0(r)\cap \M_F(1)\}\}$ and $t_i\geq 0$. Now it suffices to check that $A$ can be expressed as a convex combination of matrices in $2\{\text{conv}\{\M_0(r)\cap \M_F(1)\}\}$, namely $A=\sum_{i\geq 1}t_iB_i$. Since the zero matrix contains in $2\{\M_0(r)\cap \M_F(1)\}$, it suffices to show that $\sum_{i\geq 1}t_i\leq 1$, which is equivalent to $\sum_{i\geq 1}\|\text{vec}(D)_{T_i}\|_2\leq 2$. To prove this, first note that $\|\text{vec}(D)_{T_1}\|_2\leq \|\text{vec}(D)\|_2$. Second, note that for $i\geq 2$, each elements of $\text{vec}(D)_{T_i}$ is bounded in magnitude by $\|\text{vec}(D)_{T_{i-1}}\|_1/r$, and thus $\|\text{vec}(D)_{T_i}\|_2\leq \|\text{vec}(D)_{T_{i-1}}\|_1/\sqrt{r}$. Combining these two facts, one has that
\begin{equation*}
\sum_{i\geq 1}\|\text{vec}(D)_{T_i}\|_2\leq 1+\sum_{i\geq 2}\|\text{vec}(D)_{T_i}\|_2\leq 1+\sum_{i\geq 2}\|\text{vec}(D)_{T_{i-1}}\|_1/\sqrt{r}\leq 1+\|\text{vec}(D)\|_1/\sqrt{r}\leq 2.
\end{equation*}
The proof is complete.
\end{proof}

\begin{Lemma}\label{lem-contr}
Let $r\geq 1$, $\delta>0$ be a tolerance, and $\Gamma\in \R^{d_1\times d_1}$ be a fixed matrix. Suppose that the following condition holds
\begin{equation}\label{eq-contr1}
|\inm{\Gamma\Delta}{\Delta}|\leq \delta,\quad \forall\Delta\in \K(2r).
\end{equation}
Then we have that
\begin{equation}\label{eq-contr2}
|\inm{\Gamma\Delta}{\Delta}|\leq 12\delta(\normm{\Delta}_F^2+\frac{1}{r}\normm{\Delta}_*^2),\quad \forall\Delta\in \R^{d_1\times d_2}.
\end{equation}
\end{Lemma}
\begin{proof}
We begin with establishing the inequalities
\begin{subequations}\label{eq-contr3}
\begin{align}\
|\inm{\Gamma\Delta}{\Delta}|\leq 12\delta\normm{\Delta}_F^2,\quad \forall\Delta\in \C(r),\label{eq-contr31}\\
|\inm{\Gamma\Delta}{\Delta}|\leq \frac{12\delta}{r}\normm{\Delta}_*^2,\quad \forall\Delta\notin \C(r),\label{eq-contr32}
\end{align}
\end{subequations}
then \eqref{eq-contr2} then follows directly.

Now we turn to prove \eqref{eq-contr3}. By rescaling, \eqref{eq-contr31} holds if one can check that
\begin{equation}\label{eq-contr4}
|\inm{\Gamma\Delta}{\Delta}|\leq 12\delta,\quad \text{for all}\ \Delta\ \text{satisfying}\ \normm{\Delta}_F=1\ \text{and}\ \normm{\Delta}_*\leq \sqrt{r}.
\end{equation}
It then follows from Lemma \ref{lem-inter} and continuity that \eqref{eq-contr4} can be reduced to the problem of proving that
\begin{equation*}
|\inm{\Gamma\Delta}{\Delta}|\leq 12\delta,\quad \forall \Delta\in 2\text{conv}\{\K(r)\}=\text{conv}\{\M_0(r)\cap \M_F(2)\}.
\end{equation*}
For this purpose, consider a weighted linear combination of the form $\Delta=\sum_it_i\Delta_i$, with weights $t_i\geq 0$ such that $\sum_it_i=1$, $\normm{\Delta_i}_0\leq r$, and $\normm{\Delta_i}_F\leq 2$ for each $i$. Then one has that
\begin{equation*}
\inm{\Gamma\Delta}{\Delta}=\inm{\Gamma(\sum_it_i\Delta_i)}{(\sum_it_i\Delta_i)}=\sum_{i,j}t_it_j\inm{\Gamma\Delta_i}{\Delta_j}.
\end{equation*}
On the other hand, it holds that for all $i,j$
\begin{equation}\label{eq-contr6}
|\inm{\Gamma\Delta_i}{\Delta_j}|=\frac{1}{2}|\inm{\Gamma(\Delta_i+\Delta_j)}{(\Delta_i+\Delta_j)}-\inm{\Gamma\Delta_i}{\Delta_i}-\inm{\Gamma\Delta_j}{\Delta_j}|.
\end{equation}
Noting that $\frac{1}{2}\Delta_i$, $\frac{1}{2}\Delta_j$, $\frac{1}{4}(\Delta_i+\Delta_j)$ all belong to $\K(2r)$, and then combining \eqref{eq-contr6} with \eqref{eq-contr1}, we have that
\begin{equation*}
|\inm{\Gamma\Delta_i}{\Delta_j}|\leq \frac{1}{2}(16\delta+4\delta+4\delta)=12\delta,
\end{equation*}
for all $i,j$, and thus $\inm{\Gamma\Delta}{\Delta}\leq \sum_{i,j}t_it_j(12\delta)=12\delta(\sum_it_i)^2=12\delta$, which establishes \eqref{eq-contr31}.
As for inequality \eqref{eq-contr32}, note that for $\Delta\notin \C(r)$, one has that
\begin{equation}\label{eq-contr7}
\frac{|\inm{\Gamma\Delta}{\Delta}|}{\normm{\Delta}_*^2}\leq \frac{1}{r}\sup_{\normm{U}_*\leq \sqrt{r},\normm{U}_F\leq 1}|\inm{\Gamma U}{U}|\leq \frac{12\delta}{r},
\end{equation}
where the first inequality follows by the substitution $U=\sqrt{r}\frac{\Delta}{\normm{\Delta}_*}$, and the second inequality is due to the same argument used to prove \eqref{eq-contr31} as $U\in \C(r)$. Rearranging \eqref{eq-contr7} yields \eqref{eq-contr32}. The proof is complete.
\end{proof}

\begin{Lemma}\label{lem-rs}
Let $r\geq 1$ be any constant. Suppose that $\hat{\Gamma}$ is an estimator of $\Sigma_x$ satisfying
\begin{equation*}
|\inm{(\hat{\Gamma}-\Sigma_x)\Delta}{\Delta}|\leq \frac{\lambda_{\min}(\Sigma_x)}{24},\quad \forall \Delta\in \K(2r).
\end{equation*}
Then we have that
\begin{equation*}
\inm{\hat{\Gamma}\Delta}{\Delta}\geq \frac{\lambda_{\min}(\Sigma_x)}{2}\normm{\Delta}_F^2-\frac{\lambda_{\min}(\Sigma_x)}{2r}\normm{\Delta}_*^2,
\end{equation*}
\begin{equation*}
\inm{\hat{\Gamma}\Delta}{\Delta}\leq \frac{3\lambda_{\max}(\Sigma_x)}{2}\normm{\Delta}_F^2+\frac{\lambda_{\min}(\Sigma_x)}{2r}\normm{\Delta}_*^2.
\end{equation*}
\end{Lemma}
\begin{proof}
Set $\Gamma=\hat{\Gamma}-\Sigma_x$ and $\delta=\frac{\lambda_{\min}(\Sigma_x)}{24}$. Then Lemma \ref{lem-contr} is applicable to concluding that
\begin{equation*}
|\inm{(\hat{\Gamma}-\Sigma_x)\Delta}{\Delta}|\leq \frac{\lambda_{\min}(\Sigma_x)}{2}(\normm{\Delta}_F^2+\frac{1}{r}\normm{\Delta}_*^2),
\end{equation*}
which implies that
\begin{equation*}
\inm{\hat{\Gamma}\Delta}{\Delta}\geq \inm{\Sigma_x\Delta}{\Delta}-\frac{\lambda_{\min}(\Sigma_x)}{2}(\normm{\Delta}_F^2+\frac{1}{r}\normm{\Delta}_*^2),
\end{equation*}
\begin{equation*}
\inm{\hat{\Gamma}\Delta}{\Delta}\leq \inm{\Sigma_x\Delta}{\Delta}+\frac{\lambda_{\min}(\Sigma_x)}{2}(\normm{\Delta}_F^2+\frac{1}{r}\normm{\Delta}_*^2).
\end{equation*}
Then the conclusion follows from the fact that $\lambda_{\min}(\Sigma_x)\normm{\Delta}_F^2\leq \inm{\Sigma_x\Delta}{\Delta}\leq \lambda_{\max}(\Sigma_x)\normm{\Delta}_F^2$. The proof is complete.
\end{proof}

\begin{Lemma}\label{lem-consen}
Let $t>0$ be any constant, and $X\in \R^{n\times d_1}$ be a zero-mean sub-Gaussian matrix with parameters $(\Sigma_x,\sigma_x^2)$. Then for any fixed matrix $\Delta\in \R^{d_1\times d_2}$, there exists a universal positive constant $c$ such that
\begin{equation*}
\Pro\left[\Big|\frac{\normm{X\Delta}_F^2}{n}-\E\left(\frac{\normm{X\Delta}_F^2}{n}\right)\Big|\geq t\right]\leq 2\exp\left(-cn\min\left(\frac{t^2}{d_2^2\sigma_x^4},\frac{t}{d_2\sigma_x^2}\right)+\log d_2\right).
\end{equation*}
\end{Lemma}
\begin{proof}
By the definition of matrix Frobenius norm, one has that
\begin{equation*}
\frac{\normm{X\Delta}_F^2}{n}-\E\left(\frac{\normm{X\Delta}_F^2}{n}\right)=\sum_{j=1}^{d_2}\left[\frac{\|X\Delta_{\cdot j}\|_2^2}{n}-\E\left(\frac{\|X\Delta_{\cdot j}\|_2^2}{n}\right)\right].
\end{equation*}
Then it follows from elementary probability theory that
\begin{equation*}
\begin{aligned}
\Pro\left[\Big|\frac{\normm{X\Delta}_F^2}{n}-\E\left(\frac{\normm{X\Delta}_F^2}{n}\right)\Big|\leq t\right]
&= \Pro\left\{\Big|\sum_{j=1}^{d_2}\left[\frac{\|X\Delta_{\cdot j}\|_2^2}{n}-\E\left(\frac{\|X\Delta_{\cdot j}\|_2^2}{n}\right)\right]\Big|\leq t\right\}\\
&\geq \Pro\left\{\bigcap_{j=1}^{d_2}\left\{\Big|\frac{\|X\Delta_{\cdot j}\|_2^2}{n}-\E\left(\frac{\|X\Delta_{\cdot j}\|_2^2}{n}\right)\Big| \leq \frac{t}{d_2}\right\}\right\}\\
&\geq
\sum_{j=1}^{d_2}\Pro\left[\Big|\frac{\|X\Delta_{\cdot j}\|_2^2}{n}-\E\left(\frac{\|X\Delta_{\cdot j}\|_2^2}{n}\right)\Big|\leq \frac{t}{d_2}\right]-(d_2-1)
\end{aligned}
\end{equation*}
On the other hand, note the assumption that $X$ is a sub-Gaussian matrix with parameters $(\Sigma_x,\sigma_x^2)$.
Then \cite[Lemma 14]{loh2012supplementaryMH} is applicable to concluding that there exists a universal positive constant $c$ such that
\begin{equation*}
\begin{aligned}
\Pro\left[\Big|\frac{\normm{X\Delta}_F^2}{n}-\E\left(\frac{\normm{X\Delta}_F^2}{n}\right)\Big|\leq t\right]
&\geq d_2\left(1-2\exp\left(-cn\min\left(\frac{t^2}{d_2^2\sigma_x^4},\frac{t}{d_2\sigma_x^2}\right)\right)\right) -(d_2-1)\\
&= 1-2\exp\left(-cn\min\left(\frac{t^2}{d_2^2\sigma_x^4},\frac{t}{d_2\sigma_x^2}\right)+\log d_2\right),
\end{aligned}
\end{equation*}
which completes the proof.
\end{proof}

\begin{Lemma}\label{lem-dis}
Let $t>0$, $r\geq 1$ be any constants, and $X\in \R^{n\times d_1}$ be a zero-mean sub-Gaussian matrix with parameters $(\Sigma_x,\sigma_x^2)$. Then there exists a universal positive constant $c$ such that
\begin{equation*}
\begin{aligned}
&\Pro\left[\sup_{\Delta\in \K(2r)}\Big|\frac{\normm{X\Delta}_F^2}{n}-\E\left(\frac{\normm{X\Delta}_F^2}{n}\right)\Big|\geq t\right]\\
&\leq 2\exp\left(-cn\min\left(\frac{t^2}{d_2^2\sigma_x^4},\frac{t}{d_2\sigma_x^2}\right)+\log d_2+2r(2\max (d_1,d_2)+\log(\min(d_1,d_2))\right)
\end{aligned}
\end{equation*}
\end{Lemma}
\begin{proof}
For an index set $J\subseteq \{1,2,\cdots,\min\{d_1,d_2\}\}$, we define the set $S_J=\{\Delta\in \R^{d_1\times d_2}\big|\normm{\Delta}_F\leq 1, \text{supp}(\sigma(\Delta))\subseteq J\}$, where $\sigma(\Delta)$ refers to the singular vector of the matrix $\Delta$. Then it is easy to see that $\K(2r)=\cup_{|J|\leq 2r}S_J$. Let $G=\{U_1,U_2,\cdots,U_m\}$ be a $1/3$-cover of $S_J$, then for every $\Delta\in S_J$, there exists some $U_i$ such that $\normm{\tilde{\Delta}}_F\leq 1/3$, where $\tilde{\Delta}=\Delta-U_i$. It then follows from \cite[Section 7.2]{alquier2020high} that one can construct $G$ with $|G|\leq 27^{4r\max (d_1,d_2)}$. Define $\Psi(\Delta_1,\Delta_2)=\inm{(\frac{X^\top X}{n}-\Sigma_x)\Delta}{\Delta}$. Then one has that
\begin{equation*}
\sup_{\Delta\in S_J}|\Psi(\Delta,\Delta)|\leq \max_i|\Psi(U_i,U_i)|+2\sup_{\Delta\in S_J}|\max_i\Psi(\tilde{\Delta},U_i)|+\sup_{\Delta\in S_J}|\Psi(\tilde{\Delta},\tilde{\Delta})|.
\end{equation*}
It then follows from the fact that $3\tilde{\Delta}\in S_J$ that
\begin{equation*}
\sup_{\Delta\in S_J}|\Psi(\Delta,\Delta)|\leq \max_i|\Psi(U_i,U_i)|+\sup_{\Delta\in S_J}(\frac{2}{3}|\Psi(\Delta,\Delta)|+\frac{1}{9}|\Psi(\Delta,\Delta)|),
\end{equation*}
and hence, $\sup_{\Delta\in S_J}|\Psi(\Delta,\Delta)|\leq \frac{9}{2}\max_i|\Psi(U_i,U_i)|$. By Lemma \ref{lem-consen} and a union bound, one has that there exists a universal positive constant $c'$ such that
\begin{equation*}
\begin{aligned}
\Pro\left[\sup_{\Delta\in S_J}\Big|\frac{\normm{X\Delta}_F^2}{n}-\E\left(\frac{\normm{X\Delta}_F^2}{n}\right)\Big|\geq t\right]\leq &27^{4r\max (d_1,d_2)}\times \\
&2\exp\left(-c'n\min\left(\frac{t^2}{d_2^2\sigma_x^4},\frac{t}{d_2\sigma_x^2}\right)+\log d_2\right).
\end{aligned}
\end{equation*}
Finally, taking a union bound over the $\frac{\min(d_1,d_2)}{\lfloor2r\rfloor}\leq (\min(d_1,d_2))^{2r}$ choices of set $J$ yields that there exists a universal positive constant $c$ such that
\begin{equation*}
\begin{aligned}
&\Pro\left[\sup_{\Delta\in \K(2r)}\Big|\frac{\normm{X\Delta}_F^2}{n}-\E\left(\frac{\normm{X\Delta}_F^2}{n}\right)\Big|\geq t\right]\\
&\leq 2\exp\left(-cn\min\left(\frac{t^2}{d_2^2\sigma_x^4},\frac{t}{d_2\sigma_x^2}\right)+\log d_2+2r(2\max (d_1,d_2)+\log(\min(d_1,d_2))\right).
\end{aligned}
\end{equation*}
The proof is complete.
\end{proof}

By virtue of the above lemmas, we are now at the stage to prove Propositions \ref{prop-add-rs}--\ref{prop-mis-devia}.

\begin{proof}[Proof of Proposition \ref{prop-add-rs}]
Set
\begin{equation}\label{eq-prop1-1}
r=\frac{1}{c'}\min\left(\frac{\lambda_{\min}^2(\Sigma_x)}{d_2^2\sigma_z^4},\frac{\lambda_{\min}(\Sigma_x)}{d_2\sigma_z^2}\right)\frac{n}{2\max (d_1,d_2)+\log(\min(d_1,d_2))},
\end{equation}
with $c'>0$ being chosen sufficiently small so that $r\geq 1$. Then noting that $\hat{\Gamma}_\text{add}-\Sigma_x=\frac{Z^\top Z}{n}-\Sigma_z$ and by Lemma \ref{lem-rs}, one sees that it suffices to show that
\begin{equation*}
\sup_{\Delta\in \K(2r)}\Big|\inm{(\frac{Z^\top Z}{n}-\Sigma_z)\Delta}{\Delta}\Big|\leq \frac{\lambda_{\min}(\Sigma_x)}{24}
\end{equation*}
holds with high probability. Let $D(r):=\sup_{\Delta\in \K(2r)}\Big|\inm{(\frac{Z^\top Z}{n}-\Sigma_z)\Delta}{\Delta}\Big|$ for simplicity.
Note that the matrix $Z$ is sub-Gaussian with parameters $(\Sigma_z,\sigma_z^2)$. Then it follows from Lemma \ref{lem-dis} that there exists a universal positive constant $c''$ such that
\begin{equation*}
\begin{aligned}
&\Pro\left[D(r)\geq \frac{\lambda_{\min}(\Sigma_x)}{24}\right]\\
&\leq 2\exp\left(-c''n\min\left(\frac{\lambda_{\min}^2(\Sigma_x)}{576d_2^2\sigma_z^4},\frac{\lambda_{\min}(\Sigma_x)}{24d_2\sigma_z^2}\right)+\log d_2+2r(2\max (d_1,d_2)+\log(\min(d_1,d_2))\right).
\end{aligned}
\end{equation*}
This inequality, together with \eqref{eq-prop1-1}, implies that there exist universal positive constants $(c_0,c_1)$ such that $\tau=c_0\tau_\text{add}$, and
\begin{equation*}
\Pro\left[D(r)\geq \frac{\lambda_{\min}(\Sigma_x)}{24}\right]\leq 2\exp\left(-c_1n\min\left(\frac{\lambda_{\min}^2(\Sigma_x)}{d_2^2\sigma_z^4},\frac{\lambda_{\min}(\Sigma_x)}{d_2\sigma_z^2}\right)+\log d_2\right),
\end{equation*}
which completes the proof.
\end{proof}

\begin{proof}[Proof of Proposition \ref{prop-add-devia}]
By the definition of $\hat{\Gamma}_{\text{add}}$ and  $\hat{\Upsilon}_{\text{add}}$(cf. \eqref{sur-add}), one has that
\begin{equation*}
\begin{aligned}
\normm{\hat{\Upsilon}_{\text{add}}-\hat{\Gamma}_{\text{add}}\Theta^*}_{\text{op}}
&= \normm{\frac{Z^\top Y}{n}-(\frac{Z^\top Z}{n}-\Sigma_w)\Theta^*}_{\text{op}}\\
&= \normm{\frac{Z^\top (X\Theta^*+\epsilon)}{n}-(\frac{Z^\top Z}{n}-\Sigma_w)\Theta^*}_{\text{op}}\\
&\leq \normm{\frac{Z^\top \epsilon}{n}}_{\text{op}}+\normm{(\Sigma_w-\frac{Z^\top W}{n})\Theta^*}_{\text{op}}\\
&\leq \normm{\frac{Z^\top \epsilon}{n}}_{\text{op}}+\left(\normm{\Sigma_w}_{\text{op}}+\normm{\frac{Z^\top W}{n}}_{\text{op}}\right)\normm{\Theta^*}_*,
\end{aligned}
\end{equation*}
where the second inequality is from the fact that $Y=X\Theta^*+\epsilon$, and the third inequality is due to the triangle inequality.
Recall the assumption that the matrices $X$, $W$ and $\epsilon$ are assumed to be with i.i.d. rows sampled from Gaussian distributions $\N(0,\Sigma_x)$, $\N(0,\sigma_w^2\I_{d_1})$ and $\N(0,\sigma_\epsilon^2\I_{d_2})$, respectively. Then one has that $\Sigma_w=\sigma_w^2\I_{d_1}$ and $\normm{\Sigma_w}_\text{op}=\sigma_w$.
It follows from \cite[Lemma 3]{negahban2011estimation} that there exist universal positive constant $(c_3,c_4,c_5)$ such that
\begin{equation*}
\normm{\hat{\Upsilon}_{\text{add}}-\hat{\Gamma}_{\text{add}}\Theta^*}_{\text{op}}
\leq c_3\sigma_\epsilon\sqrt{\lambda_{\max}(\Sigma_z)}\sqrt{\frac{d_1+d_2}{n}}+(\sigma_w+c_3\sigma_w\sqrt{\lambda_{\max}(\Sigma_z)}\sqrt{\frac{2d_1}{n}})\normm{\Theta^*}_*,
\end{equation*}
with probability at least $1-c_4\exp(-c_5\log (\max(d_1,d_2))$. Recall that the nuclear norm of $\Theta^*$ is assumed to be bounded as $\normm{\Theta^*}_*\leq \omega$. Then up to constant factors, we conclude that there exists universal positive constants $(c_0,c_1,c_2)$ such that
\begin{equation*}
\normm{\hat{\Upsilon}_{\text{add}}-\hat{\Gamma}_{\text{add}}\Theta^*}_{\text{op}}
\leq c_0\phi_\text{add}\sqrt{\frac{\max(d_1,d_2)}{n}},
\end{equation*}
with probability at least $1-c_1\exp(-c_2\log (\max(d_1,d_2))$.
The proof is complete.
\end{proof}

\begin{proof}[proof of Proposition \ref{prop-mis-rs}]
This proof is similar to that of Proposition \ref{prop-add-rs} in the additive noise case.
Set $\sigma^2=\frac{\normm{\Sigma_x}_\text{op}^2}{(1-\rho)^2}$, and
\begin{equation}\label{eq-prop3-1}
r=\frac{1}{c'}\min\left(\frac{\lambda_{\min}^2(\Sigma_x)}{d_2^2\sigma^4},\frac{\lambda_{\min}(\Sigma_x)}{d_2\sigma^2}\right)\frac{n}{2\max(d_1,d_2)+\log(\min(d_1,d_2))},
\end{equation}
with $c'>0$ being chosen sufficiently small so that $r\geq 1$. Note that
\begin{equation*}
\hat{\Gamma}_{\text{mis}}=\frac{1}{(1-\rho)^2}\frac{Z^\top Z}{n}-\rho\cdot\text{diag}\left(\frac{1}{(1-\rho)^2}\frac{Z^\top Z}{n}\right)=\frac{Z^\top Z}{n}\oslash M,
\end{equation*}
and thus
\begin{equation*}
\hat{\Gamma}_{\text{mis}}-\Sigma_x=\frac{Z^\top Z}{n}\oslash M-\Sigma_x=(\frac{Z^\top Z}{n}-\Sigma_z)\oslash M.
\end{equation*}
By Lemma \ref{lem-rs}, one sees that it suffices to show that
\begin{equation*}
\sup_{\Delta\in \K(2r)}\Big|\inm{((\frac{Z^\top Z}{n}-\Sigma_z)\oslash M) \Delta}{\Delta}\Big|\leq \frac{\lambda_{\min}(\Sigma_x)}{24}
\end{equation*}
holds with high probability. Let $D(r):=\sup_{\Delta\in \K(2r)}\Big|\inm{((\frac{Z^\top Z}{n}-\Sigma_z)\oslash M) \Delta}{\Delta}\Big|$ for simplicity.
On the other hand, one has that
\begin{equation*}
\Big|\inm{((\frac{Z^\top Z}{n}-\Sigma_z)\oslash M) \Delta}{\Delta}\Big|\leq \frac{1}{(1-\rho)^2}\Big|\inm{(\frac{Z^\top Z}{n}-\Sigma_z)\Delta}{\Delta}\Big|
\end{equation*}
Note that the matrix $Z$ is sub-Gaussian with parameters $(\Sigma_z,\normm{\Sigma_x}_\text{op}^2)$ \cite{loh2012supplementaryMH}. Then it follows from Lemma \ref{lem-dis} that there exists a universal positive constant $c''$ such that
\begin{equation*}
\begin{aligned}
&\Pro\left[\Big|\inm{(\frac{Z^\top Z}{n}-\Sigma_z)\Delta}{\Delta}\Big|\geq (1-\rho)^2\frac{\lambda_{\min}(\Sigma_x)}{24}\right]\\
&\leq 2\exp(-c''n\min\left((1-\rho)^4\frac{\lambda_{\min}^2(\Sigma_x)}{576d_2^2\normm{\Sigma_x}_\text{op}^4},(1-\rho)^2\frac{\lambda_{\min}(\Sigma_x)}{24d_2\normm{\Sigma_x}_\text{op}^2}\right)\\
&\quad +\log d_2+2r(2\max (d_1,d_2)+\log(\min(d_1,d_2))).
\end{aligned}
\end{equation*}
This inequality, together with \eqref{eq-prop3-1}, implies that there exists universal positive constants $(c_0,c_1)$ such that $\tau=c_0\tau_\text{add}$, and
\begin{equation*}
\Pro\left[D(r)\geq \frac{\lambda_{\min}(\Sigma_x)}{24}\right]\leq 2\exp\left(-c_1n\min\left((1-\rho)^4\frac{\lambda_{\min}^2(\Sigma_x)}{d_2^2\normm{\Sigma_x}_\text{op}^4},(1-\rho)^2\frac{\lambda_{\min}(\Sigma_x)}{d_2\normm{\Sigma_x}_\text{op}^2}\right)+\log d_2\right),
\end{equation*}
which completes the proof.
\end{proof}

\begin{proof}[proof of Proposition \ref{prop-mis-devia}]
Note that the matrix $Z$ is sub-Gaussian with parameters $(\Sigma_z,\normm{\Sigma_x}_\text{op}^2)$ \cite{loh2012supplementaryMH}. The following discussion is divided into two parts. First consider the quantity $\normm{{\Upsilon}_{\text{mis}}-\Sigma_x\Theta^*}_\text{op}$.
By the definition of $\hat{\Upsilon}_{\text{mis}}$ (cf. \eqref{sur-mis}) and the fact that $Y=X\Theta^*+\epsilon$, one has that
\begin{equation*}
\begin{aligned}
\normm{{\hat{\Upsilon}}_{\text{mis}}-\Sigma_x\Theta^*}_\text{op}
&= \frac{1}{1-\rho}\normm{\frac{1}{n}Z^\top Y-(1-\rho)\Sigma_x\Theta^*}_\text{op}\\
&= \frac{1}{1-\rho}\normm{\frac{1}{n}Z^\top (X\Theta^*+\epsilon)-(1-\rho)\Sigma_x\Theta^*}_\text{op}\\
&\leq \frac{1}{1-\rho}\left(\normm{\left(\frac{1}{n}Z^\top X-(1-\rho)\Sigma_x\right)\Theta^*}_\text{op}+\normm{\frac{1}{n}Z^\top \epsilon}_\text{op}\right).
\end{aligned}
\end{equation*}
It then follows from the assumption that $\normm{\Theta^*}_*\leq \omega$ that
\begin{equation*}
\begin{aligned}
\normm{{\hat{\Upsilon}}_{\text{mis}}-\Sigma_x\Theta^*}_\text{op}
&\leq \frac{1}{1-\rho}\left[\left(\normm{\frac{1}{n}Z^\top X}_\text{op}+(1-\rho)\normm{\Sigma_x}_\text{op}\right)\normm{\Theta^*}_*+\normm{\frac{1}{n}Z^\top \epsilon}_\text{op}\right]\\
&\leq \frac{1}{1-\rho}\left[\left(\normm{\frac{1}{n}Z^\top X}_\text{op}+(1-\rho)\normm{\Sigma_x}_\text{op}\right)\omega+\normm{\frac{1}{n}Z^\top \epsilon}_\text{op}\right].
\end{aligned}
\end{equation*}
Recall the assumption that the matrices $X$, $W$ and $\epsilon$ are assumed to be with i.i.d. rows sampled from Gaussian distributions $\N(0,\Sigma_x)$, $\N(0,\sigma_w^2\I_{d_1})$ and $\N(0,\sigma_\epsilon^2\I_{d_2})$, respectively. Then it follows from \cite[Lemma 3]{negahban2011estimation} that there exist universal positive constant $(c_3,c_4,c_5)$ such that
\begin{equation}\label{eq-prop4-3}
\normm{{\hat{\Upsilon}}_{\text{mis}}-\Sigma_x\Theta^*}_\text{op}
\leq c_3\frac{\sigma_\epsilon}{1-\rho}\sqrt{\lambda_{\max}(\Sigma_z)}\sqrt{\frac{d_1+d_2}{n}}+
\left(c_3\frac{\normm{\Sigma_x}_\text{op}}{1-\rho}\sqrt{\lambda_{\max}(\Sigma_z)}\sqrt{\frac{2d_1}{n}}+\normm{\Sigma_x}_\text{op}\right)\omega
\end{equation}
with probability at least $1-c_4\exp(-c_5\log (\max(d_1,d_2))$.
Now let us consider the quantity $\normm{({\Gamma}_{\text{mis}}-\Sigma_x)\Theta^*}_\text{op}$.
By the definition of $\hat{\Gamma}_{\text{mis}}$ (cf. \eqref{sur-mis}), one has that
\begin{equation*}
\begin{aligned}
\normm{(\hat{\Gamma}_{\text{mis}}-\Sigma_x)\Theta^*}_\text{op}
&= \normm{((\frac{Z^\top Z}{n}-\Sigma_z)\oslash M)\Theta^*}_{\text{op}}
\leq \frac{1}{(1-\rho)^2}\normm{(\frac{Z^\top Z}{n}-\Sigma_z)\Theta^*}_{\text{op}}\\
&\leq \frac{1}{(1-\rho)^2}\left(\normm{\frac{Z^\top Z}{n}}_{\text{op}}+\normm{\Sigma_z}_{\text{op}}\right)\omega
\end{aligned}
\end{equation*}
This inequality, together with \cite[Lemma 3]{negahban2011estimation}, implies that there exists universal positive constants $(c_6,c_7,c_8)$ such that
\begin{equation}\label{eq-prop4-4}
\normm{(\hat{\Gamma}_{\text{mis}}-\Sigma_x)\Theta^*}_\text{op}\leq c_6\frac{1}{(1-\rho)^2}\normm{\Sigma_x}_\text{op}\sqrt{\lambda_{\max}(\Sigma_z)}\sqrt{\frac{2d_1}{n}}\omega+\frac{1}{(1-\rho)^2}\normm{\Sigma_x}_\text{op}\omega
\end{equation}
with probability at least $1-c_7\exp(-c_8\log (\max(d_1,d_2))$.
Combining \eqref{eq-prop4-3} and \eqref{eq-prop4-4}, up to constant factors, we conclude that there exist universal positive constant $(c_0,c_1,c_2)$ such that
\begin{equation*}
\normm{\hat{\Upsilon}_{\text{mis}}-\hat{\Gamma}_{\text{mis}}\Theta^*}_{\text{op}}
\leq c_0\frac{\lambda_{\max}(\Sigma_z)}{1-\rho}\left(\frac{\omega}{1-\rho}\normm{\Sigma_x}_\text{op}+\sigma_\epsilon\right)\sqrt{\frac{\max( d_1,d_2)}{n}},
\end{equation*}
with probability at least $1-c_1\exp(-c_2\max( d_1,d_2))$.
The proof is complete.
\end{proof}


%
%

\bibliographystyle{spmpsci}      
\bibliography{sub-JOGO-220909.bbl}   


\end{document}